\documentclass{amsart}
\usepackage{amsmath, amsthm}
\usepackage{amscd}
\usepackage{amssymb}
\usepackage{amsfonts}
\usepackage{mathrsfs}
\usepackage{latexsym}
\usepackage{graphicx}
\usepackage{wrapfig}
\usepackage[all]{xy}
\usepackage{nnfootnote}

\newtheorem{prop}{Proposition}

\newtheorem{defi}{Definition}
\newtheorem{thm}{Theorem}

\newtheorem{conj}{Conjecture}
\newtheorem{remark}{Remark}
\newtheorem{question}{Question}

\newcommand{\pants}
{\raisebox{-0.05in}
{\includegraphics
[scale=1]{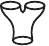}}}

\newcommand{\copants}
{\raisebox{-0.05in}
{\includegraphics
[scale=1]{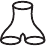}}}

\newcommand{\hH}{\mathbb{H}}

\newcommand{\lL}{\mathbb{L}}

\newcommand{\qQ}{\mathbb{Q}}
\newcommand{\rR}{\mathbb{R}}

\newcommand{\del}{\partial}

\title{String Topology \\
Background and Present State\\}
\author{Dennis Sullivan}

\begin{document}

\maketitle
\author

\nnfoottext{To appear in \textit{Current Developments in Mathematics 2005}, International Press of Boston}

\noindent
``One imagines trying to push the input circles through levels of a harmonic function on the surface. As critical levels are passed the circles are cut and reconnected near the critical points. The Poincar\'e dual cocycle creates the possibility of positioning the surface inside the target manifold.''

\begin{abstract}
(from \cite{revisedCS}) The data of a ``2D field theory with a closed string  compactification'' is an equivariant chain level action of a cell decomposition of the union of all moduli spaces of punctured Riemann surfaces with each component compactified as a pseudomanifold with boundary. The axioms on the data are contained in the following assumptions.  It is assumed the punctures are labeled and divided into nonempty sets of inputs and outputs. The inputs are marked by a tangent direction and the outputs are weighted by nonnegative real numbers adding to unity. It is assumed the gluing of inputs to outputs lands on the pseudomanifold boundary of the cell decomposition and the entire pseudomanifold  boundary is decomposed into pieces by all such factorings. It is further assumed that the action is equivariant with respect to the toroidal action of rotating the markings. A  main result  of compactified string topology is the

\begin{thm} (closed strings) Each oriented smooth manifold has a 2D field theory with a closed string compactification on the equivariant chains of its free loop space  mod constant loops.
The sum over all surface types of the top pseudomanifold chain yields a chain $X$ satisfying the master equation  $dX + X\ast X = 0$ where $\ast$ is the sum over all gluings. This structure is well defined up to homotopy*.
\end{thm}

The genus zero parts yields an infinity Lie bialgebra on the equivariant chains of the free loop space mod constant loops. The higher genus terms provide further elements of algebraic structure*  called a ``quantum Lie bialgebra'' partially resolving the involutive identity.

There is also a compactified discussion and a Theorem 2 for open strings as the first step to a more complete theory. We note a second step for knots.

\vspace{5mm}

\noindent
*\scriptsize{See the Appendix ``Homotopy theory of the master equation'' for more explanation.}

\end{abstract}

\tableofcontents

\section{\textbf{Part I. Sketch of string topology results and proofs from \cite{CS, revisedCS}}}

\subsection{The classical intersection product, infinity structures and the loop product in homology}
One knows from classical topology the homology $H_*$ of an oriented manifold of dimension $d$ has a ring structure of degree $-d$:
$$H_i \otimes H_j \stackrel{m}{\longrightarrow} H_k$$
where $k-(i+j)=-d$.

One way to define $m$ is by taking two cycles $z_1$ and $z_2$ and intersecting them transversally. In other words, after perturbing $(z_1, z_2)$ to $(z'_1, z'_2)$, intersecting $z'_1\times z'_2$ with the diagonal $M_{12}$ in $M_1 \times M_2$ where $M_1$ and $M_2$ are two copies of $M$, $((z'_1, z'_2) \pitchfork M_{12}) \subseteq M_{12} \subseteq M_1 \times M_2$.

One also knows that this intersection operation can be extended to the chain level by using a Poincar\'e dual cocycle $U$ in a neighborhood of the diagonal by considering 
$$z_1\circ z_2 = (\textrm{projection onto diagonal})((z_1 \times z_2) \cap U).$$
In this formula $\cap$ means the cap product operation $C^i \otimes C_{k+i} \rightarrow C_k$ on the chain level. On homology this intersection product is graded commutative and associative while this chain level ``diffuse intersection'' product is infinitely chain homotopy graded commutative and associative \cite{Mandell}.

Such a structure on chains is called an $E_{\infty}$ structure and such objects have a natural homotopy theory \cite{Mandell}. They can be deformed at $\qQ$ or at a prime $p$ to give higher tensors on the homology. Besides giving back the intersection product on homology, at $\qQ$ they give rise to Massey products. At $p$ they give Massey products and the Steenrod powers. At $\qQ$ and at $p$ these $E_{\infty}$ chain level structures up to homotopy determine the entire homotopy type of simply connected or even nilpotent spaces. See \cite{Su,Quillen} for $\qQ$ and \cite{Mandell} for $p$. This is true literally for closed manifolds while for manifolds with boundary it is literally true for relative chains mod the boundary.

At $\qQ$ this result may also be stated for the at $\qQ$ equivalent notion of Lie infinity or commutative infinity structure \cite{Quillen} which can also be calculated from $\qQ$-differential forms \cite{Su}. At $p$ the $E_{\infty}$ structure over the algebraic closure of the prime field must be used \cite{Mandell}.

These theorems provide a strong motivation for studying algebraic structures at the chain level and also for treating them up to homotopy. One perspective on homotopy theory of algebraic structures is sketched in the Appendix. There are others; see \cite{Costello, Toen, Vallette}.

One knows from string topology \cite{CS} how to embed this classical intersection ring structure into a ring structure on the homology $\hH_*$ of the free loop space of $M^d$:
$$\hH_i \otimes \hH_j \stackrel{m}{\longrightarrow} \hH_k$$
where $k-(i+j)=-d$.

Starting with two cycles $Z_1$ and $Z_2$ of maps of $S^1$ into $M$, evaluate at $p \in S^1$ to get two cycles $z_1$ and $z_2$ in $M$. Perturb $Z_1$ and $Z_2$ to $Z'_1$ and $Z'_2$ so that the evaluation cycles $z'_1$ and $z'_2$ are transversal to $M_{12} \subseteq M_1 \times M_2$. Then on the locus $(z'_1 \times z'_2 \pitchfork M_{12})$ in $M_{12}$ compose the loop of $Z_1$ with the loop of $Z_2$ to get the intersection cycle of loops representing $m([Z_1] \otimes  [Z_2])$.

The evaluation $Z_i \rightarrow z_i$ determines a ring homomorphism retracting the string topology product on $\hH_*$ onto the intersection product on $H_*$, where the embedding of $H_*$ in $\hH_*$ is effected by the embedding of points of $M$ into constant loops in the free loop space of $M$. These rings have a unit if and only if $M$ is closed. In the next section we define this product on the chain level.

\subsection{The chain level generalization \cite{revisedCS} of the loop product and the BV and bracket structures of \cite{CS}}

One may also define this loop product on the chain level using the diffuse intersection product $z_1 \circ z_2$. On the support of $z_1 \circ z_2$ the marked points of $z_1$ and $z_2$ determine a point in the neighborhood of the diagonal. The corresponding loops can be composed by using a short path between these nearly equal points in $M$. This chain level product is infinitely homotopy associative and so defines an $A_{\infty}$ structure in the sense of Stasheff \cite{Stasheff}. It may be transferred to an $A_{\infty}$ structure on the homology of the free loop space. There are other infinity structures on the chain level based on this diffuse string topology construction associated to further homology structures which we will now describe.

There is a circle action on the free loop space rotating the domain of maps $S^1 \rightarrow M$ which determines a degree $+1$ operator $\Delta$ on the homology. More generally, there is the $S^1$ equivariant homology $\hH^{S^1}_*$ of the free loop space and the long exact sequence relating ordinary and equivariant homology,
$$
\xymatrix{
\dots \ar[r]^M & \hH_{i+2}  \ar[r]^E
& \hH^{S^1}_{i+2}  \ar[r]^{\cap c}
& \hH^{S^1}_i  \ar[r]^M 
&\hH_{i+1}  \ar[r]^E
& \dots
}
$$
where $M$ marks a point on an equivariant loop and $E$ erases the marked point on a nonequivariant loop and, by definition, $\Delta = M \circ E$. Thus
$$ \Delta \circ \Delta = (M \circ E) \circ (M \circ E) = M\circ (E\circ M)\circ E = 0.$$

It turns out that $\Delta$ is not a derivation of the loop product $\bullet$ but the two-variable operator $\{ , \}$ defined by
$$\{x,y\} = \Delta(x \bullet y) - (\Delta x \bullet y + (-1)^{|x|} x \bullet \Delta y)$$
is a graded derivation in each variable \cite{CS}.

Since the loop product is graded commutative and associative one can say that $(\hH^{S^1}_*, \bullet, \Delta)$ forms a Batalin Vilkovisky or BV algebra since these properties are the definition of a BV algebra \cite{Getzler}. Batalin and Vilkovisky observed these structures exist on the functionals on fields of a wide variety of theories and used them to formalize quantization algorithms \cite{BV1, BV2}.

It follows from the mentioned properties of a BV algebra that $\{a,b\}$ on the ordinary homology, satisfies the jacobi identity and that the binary operation on the equivariant homology defined by
$$[a,b] = E(Ma \bullet Mb)$$
also satisfies the jacobi identity \cite{CS}. Note from the definition $\{\;,\;\}$ has degree $-d+1$, $[\;,\;]$ has degree $-d+2$, and $\hH_i^{S^1} \stackrel{M}{\rightarrow} \hH_{i+1}$ is a Lie algebra homomorphism. For $d=2$, $[\;,\;]$ is identical to the Goldman bracket for surfaces \cite{CS}. See section 2.1.

\begin{question}
How can these results from \cite{CS} be explained and what is the general picture?
\end{question}

\subsection{String diagrams for closed strings, dessins d'enfants, the combinatorial model and the general construction for the equivariant loop space}

A string diagram is a special type of ribbon graph. A \textit{ribbon graph} is by definition a graph provided with a cyclic order on the half edges at each vertex where \textit{graph} means one-dimensional CW complex. The data of the cyclic ordering allows one to construct a \textit{jet of oriented surface} along the ribbon graph for which the ribbon graph is a deformation retract or spine. One can add disks to each boundary component of the surface jet to form a closed surface with a natural cell decomposition. To be a \textit{string diagram} means the cells can be labeled input or output so no two cells of the same label meet along an edge.

String diagrams were described in terms of permutations and partitions in \cite{CS2}. They are also the same as Grothendieck's ``dessins d'enfants'' and from that context \cite{Zapponi} admit a mysterious action of the Galois group of the algebraic closure of $\qQ$.

We can form stacks \footnote{We use the word ``stacks'' as in ordinary parlance, not as in the mathematics concept.} of string diagrams where the output boundary of one is provided with a combinatorial identification to the input boundary of the diagram just below.

Stacks of string diagrams can be given a geometric interpretation by providing the input circles with metrics. Then the output circles of one string diagram inherit a metric obtained by cutting and reconnecting the input circles. This metric is transported to the next input by the identification. The choices correspond to the parameters of cells which are labeled by the combinatorial type of the stacks of string diagrams. Continuing on down, each point in the cell is represented by a cylindrical geometric surface. See figure 1.

\begin{center}
\begin{figure}[h]
\centering
\includegraphics[scale=0.8]{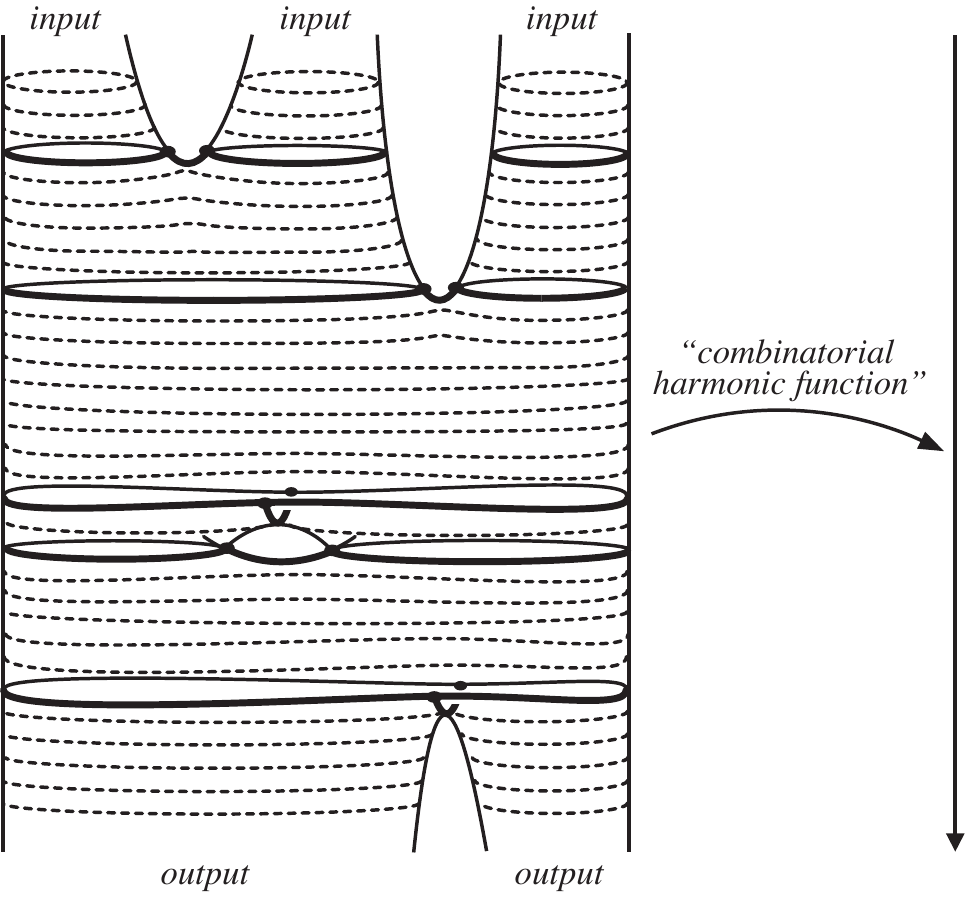}
\caption{Combinatorial harmonic function}
\label{fig:harmonicfunction}
\end{figure}
\end{center}

The proofs in \cite{CS} work with partially defined chain level operations defined using certain string diagrams that describe by transversality these operations. In \cite{revisedCS} the chain level version of \cite{CS}, narrated here, we complete the picture of this construction, replacing transversality by cap product with a local Poincar\'e dual cocycle to obtain globally defined operations on chains, and we use all possible string diagram operations.

The string diagrams when arranged in vertical stacks are separated by cylinders where heights vary in tandem between 0 and 1. The stacks exactly describe all possible compositions of cutting and reconnecting strings on the one hand with the spacings allowing deformations or homotopies between operations. On the other hand, stacks of string diagrams give a combinatorial cell decomposition of moduli spaces of Riemann surfaces as pseudomanifolds with corners which adapts to our compactification to be discussed below. 

Stacking stacks by adding a new spacing of height 1 between two stacks describes a part of the codimension one pseudomanifold boundary of these combinatorial models  of moduli space called the composition or gluing boundary.

The basic transversal string topology constructions of \cite{CS}, reformulated with Poincar\'e dual cocycles defined near the diagonal in \cite{revisedCS}, determine completely defined chain operations for each of these stacks of string diagrams (see section 3.2). These operations fit together like the cells of the model and on the composition pseudomanifold boundary give composition gluing of chain operations. Finally the models can be filled in or zipped up on the rest of the codimension one pseudomanifold boundary corresponding to small topology. We find the string topology construction extends over this extended model, in the case of closed strings.

The string topology construction starting from families of maps of the input boundary of a combinatorial surface $\Sigma$ into $M$ builds families of maps of $\Sigma$ into $M$. This reverses the direction of the left solid arrow in the following diagram. Here, $I$ and $O$ are restrictions to input and output boundary, respectively.

\vspace{5mm}
\xymatrix@C=20pt{
\{\textit{input boundary} \rightarrow M\} \ar@/^1pc/@{-->}[rr]
&
& \ar[ll]^-{I} \{\textit{surface} \rightarrow M\} \ar[rr]_-{O}
&
& \{\textit{output boundary} \rightarrow M\}
}

Then composing with the right arrow gives string-to-string operations from string-to-surface operations.

\subsection{Zipping up the noncomposition type A boundary of the combinatorial moduli space}

In all the cases, higher genus or multiple outputs or both, there will be noncomposition boundary which we want to deal with. There are two types of pseudomanifold boundary in the complement of the stacking or composition pseudomanifold boundary. For the combinatorial model we start with input circles of equal length totaling 1. As we go down through the levels and vary parameters, some part of the model carrying nontrivial topology may become small and then we reach the boundary of the model. Consider a connected subgraph $\Gamma$ of the combinatorial model surface with small combinatorial metric diameter. There are three cases depending on the Euler characteristic of $\Gamma$.  See figure \ref{fig:collapse-workhorsethm}.

\begin{enumerate}

\item[i)]
If $\chi (\Gamma) = 1$ (one less ``equation'' than ``unknowns'') we pass to a lower cell of the model, i.e., $\Gamma$ is a tree and we can collapse $\Gamma$ and stay in the moduli space of stacks of string diagrams without going beyond the pseudomanifold boundary. We are still in a lower cell of the model where the string topology construction is already defined and there is nothing further to do. (Here the approximate solutions to the ``equation'' are the support of the pulled back Poincar\'e dual cocycle  and the ``unknowns'' are the evaluations.)

\item[ii)]
 If $\chi (\Gamma) = 0$ (same number of equations as unknowns) the regular neighborhood of $\Gamma$ is homotopy equivalent to a circle and we consider the unique circuit $\alpha$ embedded in $\Gamma$ which is homotopy equivalent to $\Gamma$. The combinatorial surface may be analyzed as a local connected sum near two marked points or the circle lies on a level. Let one marked point run around a small circle with center the other in the first case; or cut along $\alpha$, twisting and regluing in the second case. This shows there is a circle action in the region of the combinatorial moduli space where $\Gamma$ has small combinatorial metric diameter, when $\chi (\Gamma) = 0$. If  $\Gamma$ is small in the combinatorial model, it is also small in $M$ because the construction (section 3.2) uses only short geodesic pieces together with pieces from the model. Thus the image of $\alpha$ may be filled by a two-disk in $M$. This means we can add a disk onto the circle orbit in the model and extend the geometric part of the string topology construction over the two-disk. The Poincar\'e part of the construction works by pull back as the combinatorial length of $\Gamma$, the radial coordinate of the two-disk, tends to zero.

This way of filling in the combinatorial moduli space and the string topology construction is familiar from the nodal curve compactification in algebraic geometry.

\item[iii)]
 If $\chi (\Gamma) < 0$ (more equations than unknowns) a convenient miracle happens. As we pull back the Poicar\'e dual cocycle to the domain of specializing evaluation -- we have, because of having more equations than unknowns, a product cocycle of bigger dimension than the specialized domain. Thus it is identically zero and the string topology construction has a vacuous locus. In this situation we can collapse to a point this part of the boundary and extend by zero. If $d>1$ this argument applies as well to chain homotopies between Poicar\'e dual cocycles. This is used in the arguments showing the ``2D field theory with closed string compactification'' is well defined up to homotopy.

\end{enumerate}

Together, these two arguments zip up the type A boundary of moduli space, where some topologically essential part of the cylindrical surface model becomes geometrically small.

\subsection{Zipping up type B boundary in the equivariant case of closed strings, indication of main theorem and its first six components}

Type B boundary: when there are multiple outputs, the input length (and the short geodesic pieces) distribute themselves among the output components. The input length distribution among the output components is described by a point in a $k$-simplex if there are $k+1$ output components. The boundary of this simplex $\times$ the combinatorial model defined without weights is by definition the type B boundary.

In the equivariant case we take care of the type B boundary by working in a quotient complex obtained by dividing out by the subcomplex of small loops. See section 3.6.

Then there is an argument that the subcomplex of equivariant chains with at least one output component very small in $M$ is mapped by any equivariant string operation into an essentially degenerate chain (section 3.6). In particular, by modding out by degenerate chains, the equivariant operations act on the quotient by the subcomplex with at least one small output component.

This removes the type B boundary from consideration in the equivariant or closed string theory.

The combinatorial model of stacks in the equivariant case has marks on the input boundary components and no marks on the output boundary components. We have weights on the outputs adding to one which record the simplex of output length distribution of the total input length. The dimension of the combinatorial model is   $-3 \chi - 1$ where $\chi$ is the Euler characteristic of the combinatorial surface $\Sigma$,  $\chi = 2-2g - \#\text{inputs} - \#\text{outputs}$.

In stacking or gluing, $\chi$ gains one by adding a mark to the output before gluing and loses one upon erasing the mark after gluing -- resulting in a net gain of zero parameters for each equivariant gluing. The fact that gluing lands in the pseudomanifold boundary of the combinatorial model is consistent with the equation 
\begin{center}
$(-3 \chi_1 - 1) + (- 3 \chi_2 -1)  = (3 \chi -1) - 1$
\end{center}
since $\chi = \chi_1 + \chi_2$.

\begin{center}
\begin{figure}[h]
\centering
\includegraphics[scale=0.75]{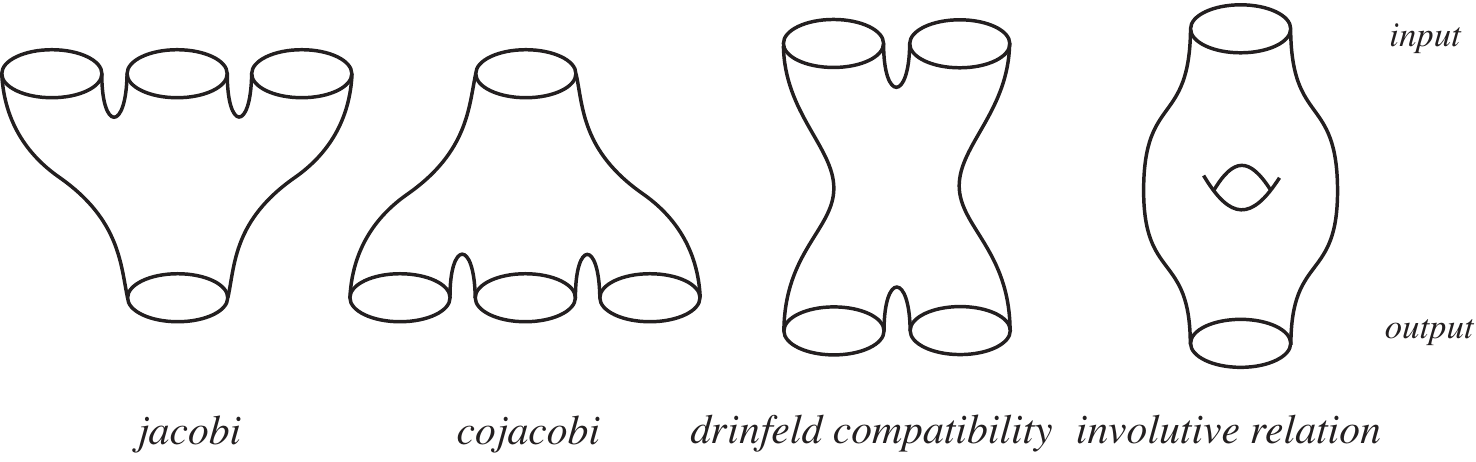}
\caption{}
\label{fig:relations}
\end{figure}
\end{center}

\noindent
\textbf{Corollary in the general equivariant case:}
In the above discussion, we have indicated there are operations on equivariant chains parametrized by the equivariant cells in a zipping up of the equivariant combinatorial moduli spaces $\mathscr{M} \subseteq \widetilde{\mathscr{M}}$ except for the composition boundary. We can apply this to the top chains of $\widetilde{\mathscr{M}}(g,k,l)$ of dimension $3 |\chi| - 1$. Let $X = \coprod_{g, k, l}$ top chain of $\widetilde{\mathscr{M}}(g,k,l)$. Then X satisfies the master equation

$$ \del X + X \ast X = 0,$$

\noindent
where $\ast$ means the sum over all possible gluings. In other words the completed chain with type A and type B boundary filled in only has composition terms in its boundary. This is the main result. It is discussed in more detail in section 3.2. We will explain in part III how the terms corresponding to \pants \hspace{0.5mm} and \copants  \hspace{0.5mm} define a bracket and cobracket  operation of degree $-d+2$ while those corresponding to  figure \ref{fig:relations}
yield four relations among these i.e.,  jacobi, cojacobi, drinfeld compatibility, and the involutive relation of an ``involutive'' Lie bialgebra.

The other surfaces give a hierarchy of higher homotopies. Restricting to the genus zero part gives a structure precisely equivalent to $\infty$ Lie bialgebra (as a dioperad). The higher genus terms give a quantum version of this infinity Lie bialgebra, called a quantum Lie bialgebra, which may be embedded in a general discussion of algebraic structures up to homotopy. (See the Appendix.)

\subsection{String diagrams for open strings, general string topology construction and the nonequivariant loop space extending \cite{CS, Su3}}

By open strings we mean oriented families of paths in an ambient manifold with endpoints on prescribed submanifolds.

For example the space of chains on the free loop space on $M$ is isomorphic to the linear space generated by equivalence classes of open strings in $M \times M$ with endpoints on the diagonal in $M \times M$.

There is a string topology construction associated to general open strings organized by the cells of a
combinatorial model of the moduli space of surfaces with boundary having marked points or punctures on the boundary. Some are designated as input open strings, some are designated as output open strings. Now we need to choose a Poincar\'e dual cocycle for each submanifold; one for the diagonal in the product of the ambient manifold and one for each diagonal in the product of each prescribed submanifold with itself. A cell that labels an operation corresponds to a combinatorial harmonic function like in figure \ref{fig:harmonicfunction} for closed strings that is proper, plus infinity at the input boundary punctures and minus infinity at the output boundary punctures. The strength of the poles at the inputs are equal to one another while they are allowed to vary at the outputs. Now there are no marked points on the input open strings.

A typical combinatorial harmonic function will have Morse-type singularities on the boundary or in the
interior and all of these are at different levels. There are two types of boundary critical points: one type increasing by two and the other type decreasing by two the number of boundary components of the level. There are two types of interior critical points: one type increasing by one, the other type decreasing by one the number of components of the level.

Starting with some input open strings, i.e., families of paths in $M$ with various boundaries, we read down the combinatorial surface using a Poicar\'e dual cocycle to impose each approximate equation at a critical point. Then as before we position the surface in the ambient manifold  by using short geodesics to build the spine.

Each critical point on the boundary imposes the number of constraints equal to the codimension of the submanifold in one case and the dimension of the submanifold in the other case. Each interior critical point imposes the ambient dimension number of constraints.By number of constraints we mean the dimension of the Poincar\'e dual cocycle that is used to impose the approximate equations.

If we glue two combinatorial surfaces together, output to input, the number of constraints is additive.

On the other hand the number of free parameters in such a top cell corresponding to a combinatorial type of combinatorial harmonic function is two for each interior critical point, one
for the increasing type of boundary critical point, zero for the the other, one for each spacing between critical levels and one less than the number of output punctures because of the variable weights. These weights must add up to the total length of the input. These quantities in total are additive less one. This dimension count checks with the fact that composition gluing puts us on the codimension one boundary of the model.

As in the closed string case, there are more types of boundary besides the composition or gluing boundary. There is the simplex boundary associated to the varying output lengths. There are also types of boundary associated to metrically small subcomplexes with nontrivial topology in the surface relative to the boundary.

The negative euler characteristic argument in the closed string case has an analogue here. The pulled-back Poincar\'e cocycle is zero for dimension reasons if the small complex has complicated topology. (A quick calculation shows we are left with collapsing circles in the interior or collapsing arcs between surface boundary as well as the simplex pseudomanifold boundary. This calculation needs to be done more carefully in general.) We treat the collapsing circles in the interior as before. Each of these remaining pseudomanifold boundary pieces of the model  can be described in terms of earlier moduli spaces.

We arrive at the result:

\begin{thm} (open strings)
The top chain in each moduli space yields an operation from strings to surfaces so that the total sum $X$  satisfies a master equation 
$$dX + X \ast X + \delta_1X + \delta_2X + \dots=0$$
where $\ast$ denotes the sum over all input output  gluings, $\delta_1$ refers to the operation inverse to erasing an output boundary puncture, $\delta_2$ refers to the operation of gluing which is inverse to cutting along the small arc,... The $\delta$ operations involve capping with Poincar\'e dual cocycles.
\end{thm}

\begin{remark}
We say more in \cite{revisedCS} about the anomalies $\delta_1$, $\delta_2 \dots$ in the nonequivariant case.
\end{remark}

\subsection{Classical knots and open string topology}

We can consider open strings in 3-space with endpoints on an embedded closed curve in 3-space.
 
With M. Sullivan we are developing an argument \cite{SullivanSullivan}to deform the $\delta$ terms to zero. This will yield a structure motivated by contact homology in the relative case. It is known from work of Ng \cite{Ng} and Ekholm, Etnyre, Ng and Sullivan \cite{EkholmEtnyreNgSullivan} that nontrivial knot invariants arise from consideration of the zeroth homology level of these invariants.

\section{\textbf{Part II. History, background, different perspectives and related work}}

\subsection{Thurston's work, Wolpert's formula, Goldman's bracket and Turaev's question}

The story of string topology begins for this author\footnote{This account, as expected, describes the story as known and remembered by the author.}  with the general background question, ``What characterizes the algebraic topology of manifolds?'' The immediate answer is the characterization should be some form of Poincar\'{e} Duality. In particular, the intersection ring of chains, 
$C_* \otimes C_* \rightarrow C_*$, defined for manifolds and its compatibility with the coalgebra structures on chains,
$C_* \rightarrow C_*\otimes C_*$, defined for all spaces, should play a role in any answer. See the CUNY theses of Mahmoud Zeinalian \cite{ZeinalianThesis} and Thomas Tradler \cite{TradlerThesis} for discussions related to duality characteristic classes and Hochschild complexes.

A second strand of the background to string topology relates to closed curves on a compact surface up to free homotopy. Their position via intersections counted geometrically rather than algebraically was important in Thurston's use of the Teichmuller space of surfaces in the study of 3D manifolds. Again, in Thurston's analysis of surface transformations he studies the orbits of embedded closed curves and how they geometrically intersect a fixed finite set of embedded closed curves instead of the usual idea in dynamics to study the orbits of points. There is also a ``cosine formula'' of Scott Wolpert for measuring the infinitesimal change of hyperbolic lengths for any geodesic $\beta$ induced by Thurston's $\alpha$ earthquake deformation, where $\alpha$ is an embedded geodesic \cite{Wolpert}. It is a weighted sum over the intersection points of $\alpha$ and $\beta$ of the cosines of the directed angles between them.

The Teichmuller space $T$ of the hyperbolic structures up to isotopy is a symplectic manifold and, by a change of variables from the cosine formula, Wolpert showed the functions $L_{\beta}$ on $T$ given by $2\cosh(\frac{1}{2} \textrm{length} \beta)$ formed a Lie subalgebra of the Poisson Lie algebra of functions on $T$ \cite{Wolpert}.

\begin{center}
\begin{figure}[h]
\centering
\includegraphics[scale=0.85]{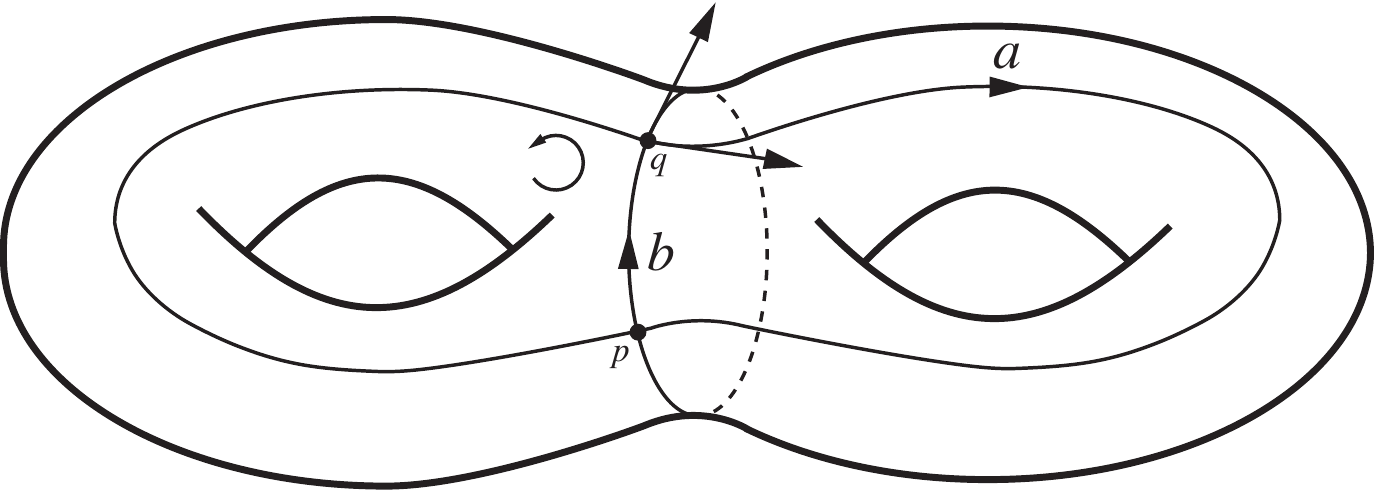}
\caption{Goldman bracket $[a,b] = (a \cdot_{q} b) - (a \cdot_{p} b)$}
\label{fig:goldman}
\end{figure}
\end{center}

Scott Wolpert suggested a Lie bracket on homotopy classes of undirected closed curves and an explanation in terms of a Lie algebra homomorphism for this remarkable Lie subalgebra fact.
 This was fully illuminated by Bill Goldman \cite{Gol} who was interested in the symplectic structure on other representation manifolds into Lie groups of the fundamental group of compact surfaces. He embedded as the invariant part under direction reversal Wolpert's construction for undirected curves into a Lie bracket
$V \otimes V \stackrel{[,]}{\rightarrow} V$ on the vector space V of directed closed curves up to free homotopy. In the example of figure \ref{fig:goldman} $[a,b] = (a \cdot_{q} b) - (a \cdot_{p} b)$ Here $(a \cdot_{x} b)$ means compose loops $a$ and $b$ at $x$ and then take the free homotopy class.

Recently, Moira Chas managed to show for embedded curves there is no cancellation in the Goldman bracket i.e., the number of terms in the bracket of two embedded closed curves is actually equal to the minimal number of intersection points \cite{Chas2}. This relates the string bracket or Goldman bracket even more forcefully to Thurston's work and suggests an algebraization is possible.

Returning again to the past, Goldman found other Lie subalgebras of functions inside the Poisson Lie algebra of all functions on the symplectic manifold of representations of the group of the surface into $G$ by considering traces, and provided an explanation using his Lie algebra of directed curves \cite{Gol}.

These manifolds of representations play a role in 3D topological quantum field theories via the geometric quantization program \cite{AndersenUeno}.

A few years after Bill Goldman's paper, Vladimir Turaev looked at self intersections of closed curves on surfaces, one by one, to split the curve into two closed curves at each self intersection. Forming a skew symmetric formal sum, Turaev defined a coLie algebra structure 
$E \rightarrow E \otimes E$ on the vector space $E$ of essential directed curves on the surface up to homotopy \cite{Turaev}. Since the trivial conjugacy class was central for Goldman's Lie algebra, the bracket of Goldman passes from 
$V$ to $E$. Note this passage from $V$ to $E$ is the precursor of the discussion above about modding out by small loops when there are multiple outputs (section 1.5).

We then  have on $E$ a Lie bracket (Goldman) 
 $E \otimes E \stackrel{[,]}{\rightarrow} E$
  and a Lie cobracket (Turaev)
 $ E \stackrel{\Delta}{\rightarrow}  E \otimes E$.
 Turaev showed the five-term drinfeld compatibility condition 
  $\Delta[a,b] = [\Delta a, b]' + [a, \Delta b]'$ where each $[,]'$, with two terms, denotes the action of the Lie algebra $E$ on its own tensor square  $E \otimes E$ \cite{Turaev}. Our generalization is sections 3.3 to 3.8 together with the new involutive property from \cite{Chas, CS2}.
  
\begin{question}
What is the deeper meaning or significance of this Lie bialgebra on the vector space $E$  of essential directed closed curves up to free homotopy? 
\end{question}

One can see the cobracket appearing in a formal expansion by Sasha Polyakov of a Wilson loop path integral calculation \cite{Polyakov}. 
  Turaev himself said he spent ten years thinking about quantizing this Lie bialgebra (which he did using knots \cite{Turaev}) and trying to understand its quantum meaning.

Turaev also asked a provocative question which led to the joint work of the author and Moira Chas ``String Topology'' \cite{CS}. Namely it  is obvious the cobracket of an embedded simple closed curve is zero and Turaev asked the beautiful question whether or not this property characterizes embedded simple closed curves among conjugacy classes which are not powers of other conjugacy classes \cite{Turaev}.

\subsection{Chas' conjectures on embedded conjugacy classes and the group theory equivalent of the Poincar\'e conjecture}

An algebraic characterization of simple conjugacy classes on 2D surfaces might be important for the topological study of 3D manifolds. There is a group theoretic statement (Jaco and Stallings \cite{Jaco, Sta}) which is equivalent to the 3D Poincar\'{e} conjecture:

``Every surjection 
$\pi_g \rightarrow F_g \times F_g$ contains in its kernel a nontrivial embedded conjugacy class.''

Here $\pi_g$ is the fundamental group of the compact surface of genus $g$, $F_g$ is the free group on $g$ generators and embedded means represented by an embedded closed curve.

 In the late '90's, motivated by this statement from the '60's, Moira Chas and the author tried to answer Turaev's question relating the kernel of the cobracket and embedded conjugacy classes. Trying to prove the affirmative answer led to the study of paths or 1-chains in the space of all closed curves on the surface. This attempt failed, but at some moment, it became clear that Goldman's bracket and Turaev's cobracket for 2D surfaces  actually existed at a geometric level in the vector spaces of chains $S_*$ of closed curves on any manifold $M^d$, 
$ S_* \otimes S_* \stackrel {[,]}{\rightarrow} S_*$ and $S_* \stackrel{\Delta}{\rightarrow} S_* \otimes S_*$ where the degree of each operation is $2-d$. Note that the degree is zero precisely for surfaces. We have indicated the full generalization above and in more detail below. See Part III.

At the same time it was clear that the basic idea, to study mapping spaces of lower dimensional manifolds into $M$ by intersecting chains in $M$ induced by evaluation at points and then forming connected sums gave a rich supply of additional operations. We chose the name ``string topology'' for the idea of intersection  followed by regluing in the study of the algebraic topology of this entire package of  mapping spaces \{intervals or circles $\rightarrow M$\}, \{surfaces $\rightarrow M$\}, etc. See e.g.,  \cite{CohenVoronov, Su4}.

We describe in more detail this ``string topology package'' for closed curves in general manifolds in Part III. Now we report on the current status of Turaev's motivating question about embedded conjugacy classes of closed curves on surfaces. 

Chas gave a combinatorial presentation of the Lie bialgebra on surfaces which could be programmed for computation and when a search was performed, examples of nonembedded  and nonpower conjugacy classes in the kernel of $\Delta$ were found (\cite{Chas}).This answered Turaev's original question about the kernel of the cobracket in the negative. 
Chas then reformulated a new conjecture about characterizing simple classes algebraically in terms of the bracket instead of the cobracket. 

\begin{conj}{(Chas)}
A conjugacy class $ \alpha $ which is not a power is simple iff any one of the following holds: 
\begin{enumerate}
\item
$[ \alpha^n, \alpha^m] = 0$ for all $n, m$
\item
$[ \alpha^n, \alpha^m] = 0$ for some $n \neq m$ and $nm \neq 0$
\end{enumerate}
\end{conj}

See \cite{Chas} for the case $n=1,m=-1$ and \cite{CK} for $n=2, m=3$.

Recently the second criterion was proven for $n=2, m=3$ by Chas and Krongold \cite{CK} so the first characterization is also established. They also suggest that a replacement for Turaev's condition $cobracket(x) = 0$, namely $cobracket(x^2) = 0$, may be sufficient to characterize embedded conjugacy classes.

Now we have Perelman's Ricci flow completion of Hamilton's program verifying Thurston's geometrization picture of 3D manifolds. Since Thurston's picture includes the Poincar\'{e} conjecture, we know that the group theoretic statement of Jaco and Stallings about embedded conjugacy classes is actually true! Thus the above-mentioned Chas-Turaev characterizations and conditions are relevant in a new way: one wants to find a purely group theory and/or topological proof of a known statement about groups which is now proved using hard PDE and hard geometry.

\subsection{Algebra perspective on string topology}

This definition of the string bracket \cite{CS} was so direct we thought it must be already known in some form. One idea possibly lay in algebra. Hochschild, in the '40's , following the idea of Eilenberg and Maclane ('40's) that groups $\Gamma$ had homology or cohomology with coefficients in any $\Gamma$-module, showed that associative algebras $A$ had homology or cohomology with coefficients in any $A$-bimodule \cite{Hochschild}.

These were defined as for groups using free resolutions and Hochschild gave specific resolutions yielding the now so-called Hochschild complexes. During the '80's it was learned \cite{Jones, GetzlerJonesPetrack} that taking the algebra to be the cochains $C^*$ on a simply connected space $X$ and the bimodule to be the chains $C_*$, the Hochschild cochains gave a model for the chains on the free loop space of $X$.

In the '80's Connes introduced the fruitful cyclic subcomplex of this particular Hochschild complex with its extra cyclic symmetry and defined the cyclic cohomology of an associative algebra \cite{Connes, Loday}. He further related this cyclic symmetry to spaces with $S^1$ action.

This cyclic structure in the Hochschild complex of $C^*(X)$ fit with the $S^1$ action on the free loop space so the cyclic cohomology of $C^*(X)$ became the equivariant homology of the free loop space of $X$ \cite{Jones}, again in the simply connected case. 

This is only half of the relevant Hochschild story. The rest occurs in the Hochschild complex for the other obvious bimodule (the predual) studied by Gerstenhaber in the '60's. 

In a celebrated paper \cite{Gerst} Gerstenhaber, motivated by a major theory at the time, the Kodaira-Spencer theory of deformations of complex structures, tried for a purely algebraic formulation. Gerstenhaber studied formal deformations of the multiplication in an associative algebra, bearing in mind the complex structure resides  in the algebra structure of the sheaf of holomorphic functions. He made use of the Hochschild complex of an algebra $A$ with coefficients in the bimodule $A$ itself. Gerstenhaber defined a $\ast$-operation in the Hochschild cochain complex $(A;A)$ which was non associative but whose commutator gave a differential Lie algebra structure on this Hochschild complex.

Nowadays one says this Gerstenhaber differential Lie algebra controls the deformation theory of $A$ and one says the Kodaira-Spencer differential Lie algebra of $(\ast, 0)$ polyvector fields with coefficients in the $(0,\ast)$ forms controls the deformation theory of the complex structure. The obstruction in each theory to extending a linear deformation $\alpha$ (here $d\alpha$ = 0 and $\alpha$ is taken mod boundaries) to a second order deformation is $[\alpha, \alpha]$ mod boundaries. If this first obstruction vanishes $[\alpha,\alpha] = d\beta$, the next obstruction to a third order deformation is $[\alpha, \beta]$ (which is a cycle in characteristic $\neq 3$ by the jacobi identity) mod boundaries, etc. These have the same universal form in any deformation theory controlled by a differential Lie algebra. Gerstenhaber also showed the Hochschild complex $(A;A)$ had a rich supply of other operations: an associative product, brace operations extending $\ast$ and so on. 

This entire Gerstenhaber discussion can be applied to the cochains $C^*(X)$ on a space. I don't know a topological interpretation of all this Gerstenhaber structure except when the space is a manifold, where it is part of the string topology of the free loop space of the manifold. 

This happens because of Poincar\'e duality at the level of chains and cochains: \textbf{in the case of a manifold $X$ the two bimodules over the cochain algebra, chains and cochains, are equivalent}. Equivalent means their resolutions are chain homotopy equivalent as chain complexes of bimodules over the algebra of cochains. This concept and fact for manifolds appeared in the CUNY thesis of Thomas Tradler \cite{TradlerThesis} which began the algebraic interpretation of string topology. Also there is an appropriate formulation for manifolds with boundary.

Thus, at least for a manifold $M$, the entire package of Gerstenhaber operations in his deformation theory may be translated from $Hochschild(C^*(X); C^*(X))$ to $Hochschild(C^*(X); C_*(X))$. The latter maps into the chains on the free loop space where the corresponding operations can be defined geometrically by the constructions of string topology. In the simply connected case the map is a quasi-isomorphism and we have an algebraic interpretation of at least a part of string topology in terms of the two Hochschild complexes $Hochschild(A;A)$ and $Hochschild(A;A^{dual})$ and the duality equivalence between them $A \backsim A^{dual}$ as $A$-bimodules in the homotopy category of $A$-bimodules. Such an algebra with this equivalence may be called a \textit{homotopy Frobenius algebra}.

Let us go back a bit before going on because we skipped a step in the story. The Lie bracket of Gerstenhaber in the above discussion is really defined on the Hochschild cohomology which models the ordinary homology $\hH_*$ of the free loop space and so defines a Lie bracket there which turns out to have degree $-d+1$. The generalized string bracket above generalizing the Lie bracket on surfaces begun by Wolpert and Goldman and generalized in \cite{CS} is defined on the equivariant homology $\hH^{S^1}_*$ of the free loop space and has degree $-d+2$. 

There is, as mentioned above, a map of degree $+1$, $\hH^{S^1}_* \stackrel{M}{\longrightarrow} \hH_*$ which is part of the long exact sequence
 $$
\xymatrix{
\dots \ar[r]^M & \hH_{i+2}  \ar[r]^E
& \hH^{S^1}_{i+2}  \ar[r]^{\cap c}
& \hH^{S^1}_i  \ar[r]^M 
&\hH_{i+1}  \ar[r]^E
& \dots
}
$$
\noindent
relating ordinary homology and equivariant homology of any space with a circle action. The map $M$ is a map of graded Lie algebras, providing a connection of the geometrically defined string bracket with the Gersenhaber bracket from algebra. 

This is a complicated path. It turns out there is a more direct route to the string bracket generalizing Goldman's bracket in surfaces. We can translate Gerstenhaber's associative product on $Hochschild(A;A)$ into a geometric product on the ordinary chains of the free loop space -- simply intersect an $\tilde{a}-$family of marked strings with a $\tilde{b}-$family of marked strings at the marked points to obtain a locus $\tilde{c}$ of dimension $c=a+b-d$. Along $\tilde{c}$ compose the loops to construct  an ``associative'' product on the chains of the free loop space of any manifold. This product $\bullet$, which we call the loop product (the Chas-Sullivan product in the literature), can be used to reinterpret the geometric definition of the string bracket using the chain maps for $E$ and $M$ of the exact sequence relating ordinary and equivariant homology (above).
$$\textit{string bracket} (a, b) = E \circ \textit{loop product} (Ma , Mb).$$
Thus begins a sequel to Gerstenhaber's algebraic deformation theory when the associative algebra satisfies some kind of Poincar\'e duality like the homotopy Frobenius property above. This occurs because the rich deformation story on one Hochschild complex (Gerstenhaber) can be combined using the equivalence of bimodules $A \backsim A^{dual}$, with the similarly interesting cyclic story (Connes) on the other Hochschild complex. 

The first combined algebraic structure that emerges is a Batalin Vilkovisky algebra structure $(\bullet, \Delta)$ on the ordinary homology of the free loop space \cite{CS} or its model the Hochschild cohomology of cochains $C^*$ with coefficients in $C^*$ or $C_*$. Here $\bullet$ is the loop product and $\Delta$ is the generator of the circle action which satisfies $\Delta \circ \Delta = 0$ and is related to the product $\bullet$ by $\{a,b\}$ which satisfies jacobi and is a derivation in each variable where $\{ , \}$ is defined by
$$-\{a,b\} = (\Delta a)\bullet b \pm a \bullet (\Delta b) - \Delta(a \bullet b).$$
In other words, the odd Poisson algebra $(\bullet , \{ , \})$ (also called a Gerstenhaber algebra) that Gerstenhaber discovered on the Hochschild cohomology $(A;A)$ for an associative algebra, in the presence of Poincar\'e duality or appropriate homotopy Frobenius, is actually derived from a BV algebra $(\bullet, \Delta)$. This algebraic statement trying to explain the geometric version of BV in string topology \cite{CS} appeared first  to my knowledge in the CUNY thesis of Tradler \cite{TradlerThesis}. 

Stronger versions, improvements and variants have appeared in several other works: Merkulov \cite{Merkulov}, R. Kaufmann \cite{Kaufmann}, Tradler \cite{Tradler}, Tradler-Zeinalian \cite{TZ}.

Currently, there are several works in progress using ribbon graphs to explore the full implications of variants of homotopy Frobenius properties of associative algebras and their generalizations to linear categories, $A_\infty$ algebras and linear $A_\infty$ categories. Kevin Costello \cite{Costello2} and Maxim Kontsevich \cite{Kontsevich2} in parallel work (suggested by earlier work of Kontsevich \cite{Kontsevich3}) have pursued the idea that the celebrated topological string theories

\begin{enumerate}
\item $A-$model defined by $J-$holomorphic curves (Gromov-Witten)
\item $B-$model defined by generalized Kodaira-Spencer deformations of Calabi-Yau manifolds
\end{enumerate}

\noindent
are conjecturally just this full algebraic exploration of the appropriate homotopy Frobenius property applied to appropriate linear ($A_\infty$) categories

\begin{enumerate}
\item for the $A-$model,  the Fukaya $A_\infty$ category of transversal Lagrangian submanifolds with morphisms given by the Floer complexes of $J-$holomorphic curves
\item for the $B-$model, the $A_\infty$ category of coherent sheaves on a Calabi-Yau manifold where the morphisms are the complexes for defining $Ext( , )$.
\end{enumerate}

In the examples $A-$model and $B-$model it is important that the cyclic Hochschild complex here has finite type and satisfies Poincar\'e duality in its own right which does not occur for string topology because the free loop space is infinite dimensional.

V\'eronique Godin \cite{Godin} and Ralph Kaufmann \cite{Kaufmann2, Kaufmann3} and Mike Hopkins with Jacob Lurie \cite{HopkinsLurie} are developing this ribbon graph picture of what can be called the ``$B_\Gamma$ part'' of string topology from either the original string topology perspective or this homotopy Frobenius algebra perspective. 

By the ``$B_\Gamma$ part'' of string topology we mean the part of the structure described below associated to the cells in the interior of the moduli space of Riemann surfaces together with those on the boundary associated to gluing Riemann surfaces. 

By $\Gamma$ we mean the mapping class group and one knows (in several different statements) that the open part of moduli space is homologically equivalent to the classifying space $B_\Gamma$ of $\Gamma$. 

As we mentioned above, there is more to discuss about the $B_\Gamma$ part of string topology. The issue is whether or not the rest of the boundary of the open moduli space beyond composition can be ``zipped up'' or compactified. 

In compactified string topology \cite{revisedCS} the boundary is zipped up by combinatorially finding and then filling in loops that are small in the manifold.

\subsection{Homotopy theory or algebraic topology perspective on string topology}

\textbf{Umkehr map of string topology} The operations in string topology use a wrong way or umkehr map associated to the left arrow in the diagram

\xymatrix@C=20pt{
(\textit{input boundary}, M)
&
& \ar[ll]^-{I} (\textit{surface}, M) \ar[rr]_-{O}
&
& (\textit{output boundary}, M)
}

Here $(X,M)$ is the space of all maps of $X$ into $M$. This ``umkehr map'' is defined after applying a linearizing functor, say $F$, to the diagram and then doing a version of intersection product to get the umkehr map

\xymatrix{
F(\textit{input boundary},M) \ar[rr]^-{umkehr}
& &
F(\textit{surface},M).
}

For example, taking the surface to be a pair of pants with marked input boundary (2 in, 1 out) and $F$ to be the chains on the loop space, leads to the loop product. If $F$ is the subcomplex of equivariant chains inside all the chains, one gets the string bracket. Here umkehr is geometrically intersecting with the diagonal or, as we mentioned above and described in more detail in sections 1.3 and 3.1 to 3.8, pulling back a Poincar\'e dual cocycle to the diagonal, e.g., a Thom class representative for the (neighborhood of the diagonal, boundary). 

One could also take $F$ to be representatives of bordism instead of chains. The most general object to use instead of the chains is the spectrum of any cohomology theory $h^*$ for which the normal bundle of the diagonal $M \rightarrow M \times M$ is orientable. 

Cohen and Jones \cite{CJ} have devised such a stable homotopy formulation of some of the operations in the $B_\Gamma$ part of string topology corresponding to string diagrams. Recently, V\'eronique Godin has developed a spectrum formulation of the ``$B_\Gamma$ part'' of string topology \cite{Godin} (see previous section for a definition of ``$B_\Gamma$ part''). Jacob Lurie and Mike Hopkins have a categorical understanding of a general form of these constructions and a potential framework for the $A$ and $B$ models mentioned in the previous section \cite{HopkinsLurie}. I suppose their framework will eventually include the compactified form of string topology discussed here. 

These constructions respect the composition part of the boundary of open moduli space. I don't know what to expect about their behavior near the rest of the boundary of open moduli space at infinty, except to say that the Euler class, which is the image of the class of the Poincar\'e dual cocycle under the map
$$h^d(\textit{neighborhood of diagonal}, \del(\textit{neighborhood of diagonal})) \rightarrow h^d(\textit{neighborhood of diagonal})$$
\noindent
should play a role. 

The above remarks concerned stable homotopy theory. Here is  a connection to unstable homotopy theory.

The $E_2$ term of the homology spectral sequence for the natural fibration
$$
\xymatrix{
\textit{based loop space} \ar[rr]^{inclusion} && \textit{free loop space} \ar[rr]^{evaluation} && \textit{manifold}
}
$$
\noindent
is the tensor product of the homology of the based loop space and the homology of the manifold ($\qQ$ coefficients, no twisting). The first factor is a graded cocommutative Hopf algebra and the second is a graded commutative Frobenius algebra. 

The tensor product $E_2$ term is therefore both an algebra and a coalgebra. It is not clear how to express an intelligent compatibility between these two structures. 

Now the differentials of the spectral sequence respect the coalgebra structure (true for all spaces) and this leads to the coalgebra structure on the limit which agrees with the coalgebra structure on the free loop space. Before string topology came along one could have asked if the differential $d_3$ preserves the algebra structure on $E_2$ and then ask if this algebra structure persists through the spectral sequence to some algebra structure on the homology of the free loop space. 

Actually, the geometric construction of the loop product at the chain level respects the filtrations defining the spectral sequence. Thus the loop product exists all through the spectral sequence and the differentials respect both the algebra and the coalgebra structure. 

Cohen, Jones and Yan \cite{CohenJonesYan} also noticed and then emphasized this neat fact and used this property to obtain simple computations of the free loop space homology for familiar spaces. 

Recently, Xiaojun Chen, in his Stony Brook thesis, \cite{ChenThesis} has built a chain model of this free loop space fibration ($\qQ$ coefficients) with an explicit differential on the tensor product of the forms and a completed cobar construction on the dual of forms that is both a derivation and coderivation for the product and completed coproduct. 

An intriguing problem is to formulate the kind of bialgebra this construction instantiates. Up to now the diagonal coalgebra structure of the free loop space has stood somewhat apart from the algebraic structures on the free loop space coming from string topology. A special case of the problem is illustrated by the $E_2$ term above -- how should one view the tensor product of a Frobenius algebra and a Hopf algebra? 

Of course the Hopf algebra there is really the universal enveloping algebra of a Lie algebra. The tensor product of a commutiative algebra with a Lie algebra is a Lie algebra. Also our commutative algebra has an invariant trace in the closed manifold case. But then what? 

\begin{question}
How much of string topology is a homotopy invariant of the pair $(M, \del M)$? In \cite{CKS} it is shown that the string bracket and the loop product are homotopy invariants of closed manifolds. Perhaps the entire $B_\Gamma$ part can be constructed from the homotopy theory, using Hopkins' and Lurie's construction \cite{HopkinsLurie}. 
\end{question}

On the other hand we have conjectured that the entire ``zipped up'' string topology package is not a homotopy invariant (see Postscript \cite{Su2} and next section). The motivation for this conjecture can be obtained from the sequence of statements.

\begin{enumerate}
\item
(Ren\'e Thom '58) A polyhedron which is locally a $\qQ-$homology manifold has $\qQ-$Pontrjagin classes \cite{Thom}. These classes are not homotopy invariants but in fact parametrize the infinite part of the diffeomorphism types for higher dimensional simply connected manifolds ('60's surgery theory).
\item
(Clint McCrory '70) An oriented pseudomanifold $P$ without boundary for which the diagonal $P \rightarrow P \times P$ has a dual cocycle supported in a regular neighborhood of the diagonal (mod its boundary) is a homology manifold \cite{McCrory}.
\item
The string topology constructions discussed above precisely use a local near the diagonal Poincar\'e dual class to construct the chain level string topology operations and the locality seems necessary.
\end{enumerate}

\subsection{Symplectic topology perspective on string topology}

There has been renewed activity (referred to as symplectic topology) since 1985 and Gromov's discovery of the control plus flexibility of $J-$holomorphic curves (i.e., surfaces mapping into a symplectic manifold provided with an almost complex structure $J$). Homological invariants of the moduli spaces of such curves with specified boundary conditions and lying in fixed relative 2-dimensional homology classes provide a rich array of invariants. Being homological with $\qQ-$coefficients, these invariants remain unchanged as the almost complex structure and the symplectic structure are deformed continuously. Thus they can in this sense be considered to be part of topology (as the physicists have done for years).

These theories in various forms can be applied to general smooth manifolds $M$ by considering the cotangent bundles and their tautological exact symplectic structures ($\omega = d\eta$ where locally $\eta = \sum_i p_i dq_i$). 

There are at least three forms of symplectic topology that may be used here:

\begin{enumerate}
\item
a Floer type theory that leads to  operations in a $J-$holomorphic disk description of the ordinary homology of the free loop space of $M$ \cite{Viterbo, Viterbo', SalamonWeber, AbbondandoloSchwarz, Cohen}. This relates to the nonequivariant loop space (open string topology).
\item
the symplectic field theory (SFT) \cite{EliashbergGiventalHofer} applied to the cotangent bundle minus the zero section regarded as the symplecticization of a contact manifold, the unit sphere cotangent bundle. $J-$holomorphic curves in the symplectization can describe the equivariant homology of the free loop space. This relates to the equivariant loop space (closed string topology).
\item
relative symplectic field theory, also related to open string topology.
\end{enumerate}

Yasha Eliashberg has emphasized two interlocking questions:

\begin{enumerate}
\item
does the symplectic structure on $T^*M$ determine $M$ up to diffeomorphism and is $Diff(M)$ homotopy equivalent to $Symplectomorphism(T^*M)$?
\item
can all the known invariants of smooth manifolds: the homotopy type, the characteristic classes, the surgery invariants of higher dimensional manifolds, the Donaldson and Seiberg-Witten invariants of 4-manifolds, and the quantum invariants of 3-manifolds of Chern-Simons, Vaughan Jones and Vasiliev, be described in terms of $J-$holomorphic curve invariants of the cotangent bundle?
\end{enumerate}

Eliashberg has also perceived a role for string topology in the general theory of $J-$holomporphic curve invariants of pairs (symplectic $W$, systems of Lagrangian submanifolds $L_1, L_2, \dots$). Namely, conjecturally, the symplectic theory of the cotangent bundle may be described in terms of string topology and also maybe this constitutes a natural piece of symplectic topology near the Langrangian boundary conditions. This might happen because of Weinstein's result that the neighborhood of any Lagrangian $L$ in symplectic a symplectic manifold $W$ is symplectomorphic to a neighborhood of the zero section in the cotangent bundle of $L$. 

Let us examine in somewhat more detail the first point of this speculation and discuss a bit the program of Janko Latschev and Kai Cieliebak \cite{CieliebakLatschev}. Starting from the algebraic formulation of SFT \cite{EliashbergGiventalHofer} they consider, in the case of the cotangent bundle, a rich algebraic invariant of $M$. In fact, there are three levels as follows. 

The formulation of SFT \cite{EliashbergGiventalHofer} uses punctured $J-$holomorphic curves in the cotangent bundle minus the zero section stretching between periodic orbits of the Reeb flow (e.g., the geodesic flow) at $+\infty$ or $-\infty$ (which is near the zero section). The genus zero curves with one component at $+\infty$ (level I) leads to a derivation differential $d=d_0+d_1+d_2+ \dots$ on the free graded commutative algebra on periodic Reeb orbits. Using the filling of the contact manifold by the unit disk bundle a change of variables can be discovered that reduces the constant term $d_0$ to zero. Then $d_1 \circ d_1=0$ and the linearized homology (the homology of $d_1$ on the indecomposables) turns out to be the equivariant homology of the free loop space mod constant loops, i.e., the reduced equivariant homology used in section 2.2. 

Analyzing the string cobracket of string topology (see below) leads to a similar structure -- a coLie infinity structure on  $\widetilde{\hH}^{S^1}_*$, the reduced equivariant homology of the free loop space. 

The Cieliebak-Latschev program at level I is to construct this type of structure by transversality as in the string topology of \cite{CS, CS2} (rather than the Poincar\'e dual cocycle version here) and show it is isomorphic to the structure coming from the $J-$holomorphic curves. 

They have a similar program for level II using genus zero curves with several punctures at $+\infty$. Now one can view the SFT formalism as an infinity Lie bialgebra structure (see our description below of the infinity Lie bialgebra structure arising from string topology). They try to construct this structure as in string topology \cite{CS, CS2} and again show it is isomorphic to the one coming from level II $J-$holomorphic curves. 

Actually the Lie bialgebra in string topology is involutive at the chain level. One now understands the infinity version of this genus one relation for a Lie bialgebra requires at least and perhaps more operations indexed by ($k$ inputs, $l$ outputs, $g=$ genus). Their level III program uses the higher genus curves as well. See section 1.5 above.

There is also relative symplectic field theory that can be applied to study classical knots $K$ in 3-space. The conormal of $K$ in the cotangent bundle of 3-space provides the boundary conditions for the $J-$holomorphic curves. The level I theory adjusted by the filling given by the relative cotangent disk bundle yields a differential derivation $d=d_1+d_2+ \dots$ on a free associative algebra generated by the Reeb flow orbits starting and ending on the Lagrangian. Lenny Ng has a set of papers motivated by the problem of computing the zeroth homology of this dga. He found a conjectured combinatorial description and showed it gives a powerful knot invariant \cite{Ng}.

This conjecture is now proved in \cite{EkholmEtnyreNgSullivan}. It turns out that Ng's combinatorial descriptions resonates with the open string topology of the knot \cite{SullivanSullivan}-- it is related to the coproduct on families of strings in $\rR^3$ starting and ending on the knot which are cut by intersecting with the knot. A Poicar\'e dual cocycle description of the intersection defining $d_2$ can be chosen to eliminate the anomaly and this should lead to an $A_\infty$ coalgebra structure on the linearlized homology. Again there are variants of the construction \cite{SullivanSullivan}.

\subsection{Riemannian geometry perspective on string topology}

One might imagine making the string topology constructions using the heat flow. Each heat operator $e^{\Delta t}$ commutes with $d$ and is chain homotopic to the identity via $\int_{0}^{t}  d^*e^{\Delta s} ds$. Also $e^{\Delta t}$ provides a differential form Poincar\'e dual to the diagonal with more and more of its weighted support tending to the diagonal as $t$ tends to zero. In the probabilistic picture this diffusion operator is related to parallel translation modified by a curvature term along random paths weighted by the Wiener measure \cite{Wiesbrock}. One could imagine using these paths to define the string topology operations instead of the short geodesics. In fact, as $t \rightarrow 0^+$ the Wiener measure conditioned on paths from  $x$ to $y$ converges to a measure on the shortest geodesics from $x$ to $y$ \cite{DeuschelStroock}. The details of this putative heat string topology construction are nontrivial but feasible (see \cite{Wiesbrock, Stasheff}). 

However, Kevin Costello \cite{Costello} has a completed diagrammatic calculation using these quantities, harmonic forms and $e^{\Delta t}$, involving ribbon graphs. Renormalization issues are addressed and a rich structure is produced. It reminds one of string topology for small loops with coefficients in a compact Lie group $G$ -- resonating with the original work of Goldman on surfaces.

\section{\textbf{Part III. The diffusion intersection and short geodesic construction of string topology, the main theorem and the first six examples}}

\subsection{Statement of the main theorem of string topology for closed strings and the motivation for infinity structures}

The generalized bracket for two families of closed strings $A$ and $B$ is very simple to define geometrically --  just intersect (transversally) the set of possible positions on all the curves in the $A$ family with the same in the $B$ family. At each point of this locus $C$, compose the $A$ curves and $B$ curves, as based loops, to define the Lie bracket family $C$ of unbased loops. A picture \cite{CS2} shows this operation satisfies jacobi and passes to a Lie bracket on the equivariant homology of the free loop space of $M^d$. The degree of this operation is $-d+2$ as can be seen from the intersection theory used above. If $A$ has dimension $a$, B has dimension $b$, then $C$ has dimension $c = a+b+2-d$. This process is exactly the formula $[a,b] = E(Ma \bullet Mb)$, mentioned in section 1.1.

One knows in algebraic topology, it is not really optimal to pass to homology in this kind of situation but one also knows that the alternative is more work. The problem and the interest is that the above bracket and jacobi identity only pertain``transversally.'' In the Stony Brook thesis of Scott Wilson \cite{WilsonThesis} such partially defined structures with one output were extended to globally defined structures on functorially associated quasiisomorphic complexes. The theory of the Appendix provides an adequate theory for multiple outputs. How this partial-to-global transition should be interpreted can be learned from the intersection product of relative integral chains in a manifold with boundary. The \textit{transversal} intersection product is \textit{actually} graded commutative and associative. Steenrod operations show it \textit{cannot} be extended to such a product on all integral chains. However, the general theory begun by Steenrod says it may be extended to a commutative and associative product up to infinite homotopy. Furthermore, as mentioned above, over $\qQ$ \cite{Quillen, Su} and over $\bar{F_p}$ for each prime $p$ \cite{Mandell} this $E_{\infty}$ product structure up to homotopy determines the entire homotopy type for simply connected spaces. In fact in each setting (at $\qQ$ or at $p$) there is an equivalence of homotopy categories between spaces and the algebraic models.

One of the main consequences of the general string topology  construction described below provides analogous Lie bracket and Lie cobracket infinity structures for the chains on the free loop space, namely one has the following theorem.

\begin{thm}
On the reduced equivariant chains of the free loop space, $\lL^{S^1}_*(k)$  of an oriented $d$-manifold (defined below in section 3.3), the ``diffuse intersection'' string topology construction produces an ``involutive Lie biagebra structure up to homotopy.'' The degrees of the bracket and cobracket operations are $2-d$.
\end{thm}

What does the theorem mean algebraically? We will explain the quotation marks in remark 2 below. There are two parts to the theorem:

\begin{enumerate}
\item[Part I]
For each $(k,l,g)$ associated to the top or fundamental chain of the combinatorial moduli space there is a well defined graded symmetric chain operation of degree 
\newline
$(3-d)(2g-2+k+l)-1$: $\phi_g: \lL^{S^1}_*(k)' \rightarrow \lL^{S^1}_*(l)', k>0, l>0, g\geq 0$. Add a formal variable $t$ in order to sum these operations obtaining: $\varphi (k,l) = \sum_{g=0}^{\infty} t^{2g-2}\phi_g$.
For $d=4,5,6, \dots$ all but finitely many of the operations $\phi_g(k,l)$ are zero in the fixed degree because their degrees go to $-\infty$.
So form $\sum t^{k+l} \varphi(k,l) = \vartheta$.
The theorem says that these operations are defined and satisfy the master equation:
$$ \del \vartheta + \vartheta \ast \vartheta = 0$$
where $\del \vartheta =$ the commmutator of $\vartheta$ with the boundary operator and $\vartheta \ast \vartheta$ is the sum over all possible compositions.

\item[Part II]
Structures like those described in Part I can be transported between different chain complexes via chain mappings inducing isomorphisms on homology (see Appendix). With $\qQ$ coefficients the homology with zero differential is such a complex. Thus from the $\vartheta$ in Part I on $C_*$ there is implied a collection \{$\vartheta_H(g,k, l)$\} acting between $(\hH^{S^1})^{\otimes k} \rightarrow (\hH^{S^1})^{\otimes l}$.

Adding the formal variable $t$ and summing as before $\vartheta_H = \sum_{k,l} t^{k+l} \sum_{g \geq 0} \varphi_g(k,l)$ to get one operation $\vartheta_H$ on $\Lambda(\hH^{S^1}_*t)$, the free graded commutative algebra on $\hH^{S^1}_*$ the reduced equivariant homology of the free loop space of $M^d$ with $t$ then adjoined. The above equation at the chain level $ \del \vartheta + \vartheta \ast \vartheta = 0$ becomes the equation $\vartheta_H \ast \vartheta_H = 0$ at the homology level in $\Lambda(\hH^{S^1}_*t)$.

If we give $t$ the weight $d-3$ then $\vartheta$ and $\vartheta_H$ each has degree $-1$.

\end{enumerate}

\begin{remark}
The algebraic structure indicated by the master equation of Part I is not a complete resolution of the involutive Lie bialgebra structure, thus the quotation marks in the theorem. It is however an infinity algebraic structure in the sense of the Appendix -- namely an infinity version of its own homotopy type. This homotopy type could be named  \textit{quantum Lie bialgebra}.
The situation is analogous to work of Kevin Costello \cite{Costello2} where the diagrams constructed give a version of a resolution of  the  cyclic $A_\infty$ structure and Costello dubs the structure a \textit{quantum $A_\infty$ structure}.
\end{remark}

\begin{remark}
For dimension $d=2$, the reduced equivariant homology of the free loop space for higher genus is concentrated in degree zero ($\qQ$ coefficients) and is just the space $E$ of the introduction. The operations $\phi(k,l,g)$ have degree $(3-d)(2g-2+k+l)-1$, which is nonzero unless $g=0$ and $k+l=3$. This leaves only the bracket and cobracket of Goldman and Turaev for $d=2$ if the surface is not the 2-sphere or the torus.

For dimension $d=3$ the degree of every operation is $-1$. For closed hyperbolic manifolds the homology is concentrated in degree zero ($\qQ$ coeffiecients) so all operations $\phi_H$ of the minimal model are zero. In the CUNY thesis of Hossein Abbaspour a converse is proven expressed in terms of the loop product on the ordinary homology of the free loop space. Namely, if a closed 3-dimensional manifold is not hyperbolic, some string topology (loop products beyond classical intersection products) is nontrivial. See the precise statement \cite{Abbaspour} where double covers must be used for certain ``small Seifert fibered spaces.''
\end{remark}

\subsection{Sketch of the basic diffuse intersection and geodesic path construction of string topology}

Start with a pair consisting of a family of oriented closed one-manifolds in $M$ with $k$ labelled components and a combinatorial description (via a combinatorial harmonic function, see \cite{CS2, Bodigheimer} and figure \ref{fig:harmonicfunction}) of a combinatorial surface of genus $g$ with $k$ labeled input $\del$ components and $l$ labeled output components. One imagines trying to push the input circles through the surface. As critical levels are passed the circles are cut and reconnected precisely at the critical points. This happens a finite number of times until the output boundary is achieved. Between critical levels one imagines only moving slightly, if at all. The diagrammatic picture of the surface changes when two critical levels come together and unite into one critical level. 

\begin{center}
\begin{figure}[h]
\centering
\includegraphics[scale=1.00]{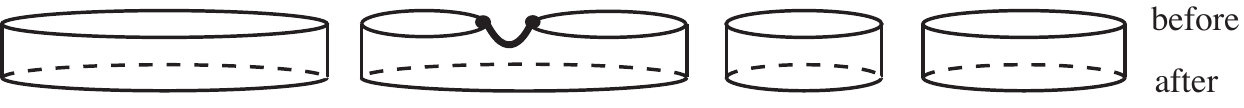}
\caption{Just before and just after a critical level in the surface}
\label{fig:criticallevelsurface}
\end{figure}
\end{center}

For the generic picture of the surface each critical level has one Morse quadratic-type critical point which is neither a minimum nor a maximum. In this case there are $(2g-2+k+l)$ critical levels since each one adds $-1$ to the Euler characteristic of what came before. The critical level and what came just before and just after is specified by two parameters: where the two points of the ascending manifold attach. See figure \ref{fig:criticallevelsurface}.

In the generic case there are $(2g-2+k+l)-1$ heights of cylinders between critical levels which we sometimes refer to as the spacing between levels. All in all there are 
\newline
$2(2g-2+k+l) + (2g-2+k+l) - 1$ parameters to describe the typical picture, so combinatorially a top cell or stratum is of dimension $3| \chi |-1$ where $\chi$ is the Euler characteristic of the punctured surface.

These cells or strata fit together to form a compact pseudomanifold with corners of that dimension. The boundary is created by imposing certain inequalities that the combinatorial metric length of any essential combinatorial circuit is $\geq \epsilon > 0$. Boundary is also created by spacing equal to one because spacing is restricted to the interval $[0,1]$.

When surfaces 1 and surface 2 are glued output to input to obtain surface 3 note that $\chi_3 = \chi_1 + \chi_2$ but the dimensions of the moduli spaces of these satisfy $d_1 + d_2 + 1 = d_3$. It turns out that the product of the 1 and 2 moduli spaces is embedded in the boundary of the 3 moduli space by gluing. This is also clear because, as mentioned above, the cylinders created by gluing have length 1 by definition.

The string topology construction will construct, on some open subset of the initial family of input boundary mapping into $M$, a mapping of the combinatorial surface into $M$ with the given input boundary values. This is done step by step, over each level. If inductively (over $r$, say) we have mapped the surface up to just below the critical level on some open set $U_r$, we map the two attaching points of the next critical level into $M \times M$.

\begin{center}
\begin{figure}[h]
\centering
\includegraphics[scale=1]{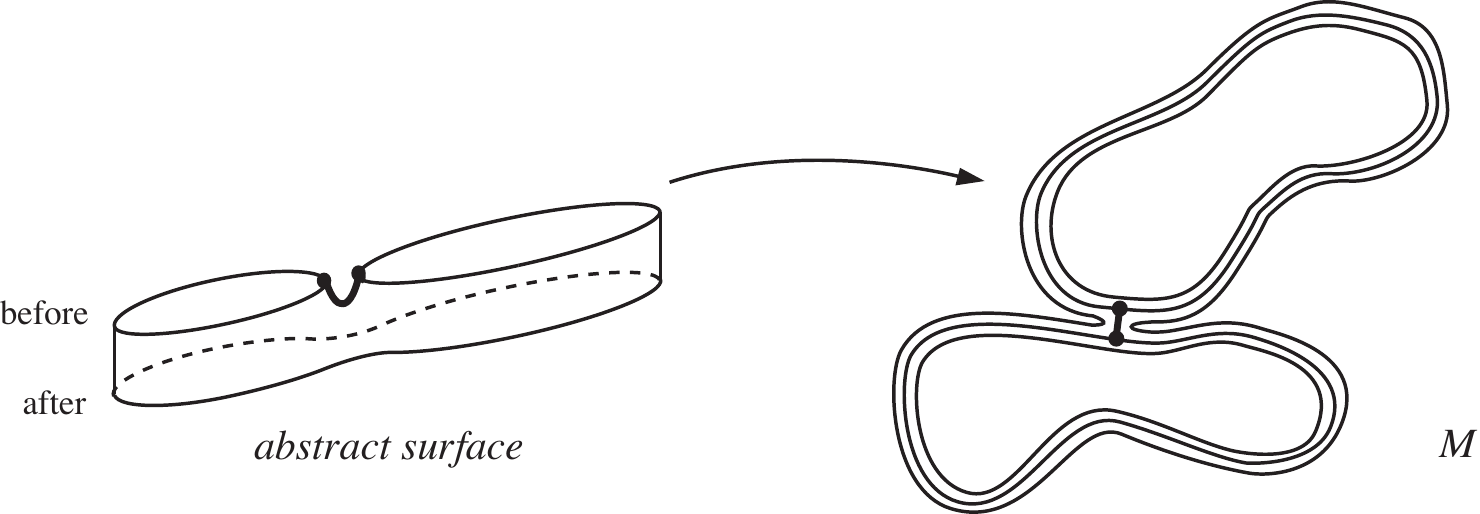}
\caption{Near a critical level in the surface and in the manifold}
\label{fig:criticallevelmanifold}
\end{figure}
\end{center}

We pull back a small neighborhood of the diagonal in $M \times M$ to define $U_{r+1}$ inside $U_r$. On $U_{r+1}$ the two points will be close in $M$. We can draw a short geodesic between these points and map the short critical trajectory onto this geodesic. The level after the critical level and just before the next critical level is projected in the surface to the level just below the critical arc plus the critical arc. Putting this together gives the map  into $M$ of the surface up to just before the next critical level. See figure \ref{fig:criticallevelmanifold}.

This describes the basic geometric construction for interiors of strata corresponding to the combinatorial harmonic function picture being Morse with distinct critical levels. We postpone for the moment how this process defines a chain operation and first discuss what happens as critical levels merge or separate. The general picture of one critical level is combinatorially equivalent to the string diagrams or generalized chord diagrams introduced in \cite{CS2} and defined above in section 1.3. The relationship of stacks of these to Riemann surfaces with a harmonic function goes back to Poincar\'e and Hilbert. (See \cite{Bodigheimer} for more details.)

\begin{center}
\begin{figure}[h]
\centering
\includegraphics[scale=1]{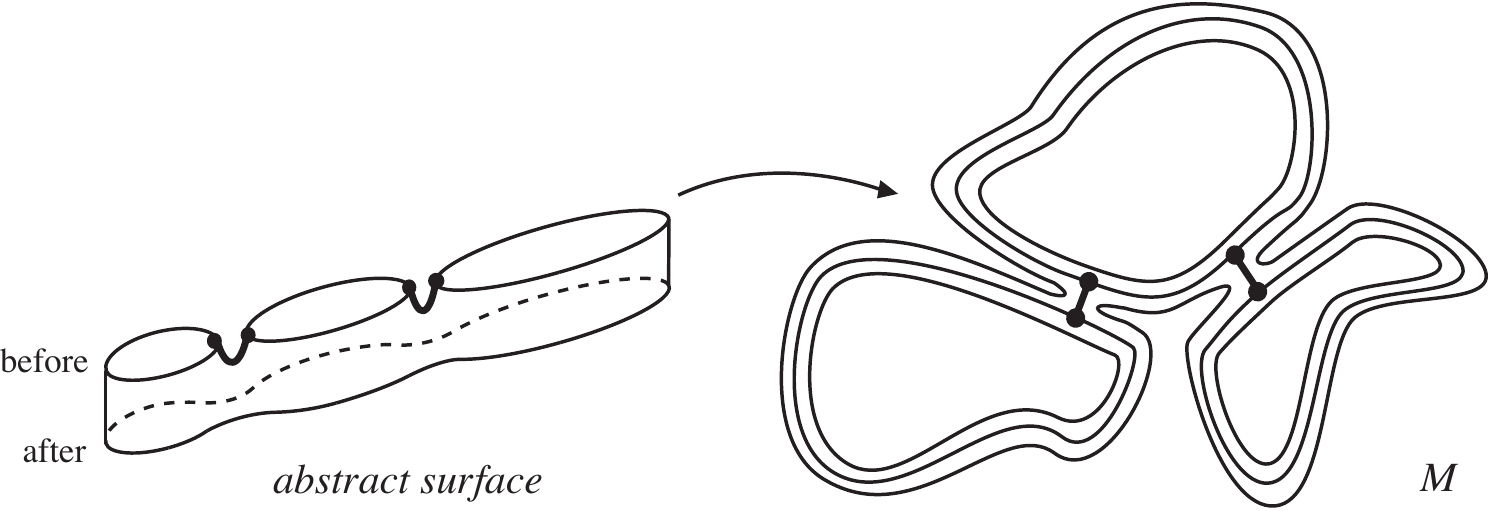}
\caption{Critical levels come together, in the surface and in the manifold}
\label{fig:criticallevelstogether}
\end{figure}
\end{center}

If two levels come together in a generic way there are distinct Morse critical points at the same level. There is an associated evaluation map to four copies of $M$ and we can pull back the intersection of two neighborhoods of the diagonals, say $12$ and $34$ in $(1,2,3,4)$. This intersection is the Cartesian product of the neighborhood of the diagonal $12$ in $(1,2)$ and the diagonal $34$ in $(3,4)$. For the geometric construction we form two short geodesics and get the picture in figure \ref{fig:criticallevelstogether}.

Note this picture is not identical to one of the other pictures obtained by doing first one operation and then the other. However these two only differ on small convex neighborhoods in $M$ of the points in question. So it is possible to construct a canonical homotopy between the two maps of the surface into $M$. We associate the parameter of this homotopy to the parameter corresponding to the height of the cylinders between critical levels.

The above geometric construction defines operations on families of input circles mapping into $M$ -- where the family is cut down to the open set where the appropriate equations approximately hold allowing us to position the surface in $M$. To get chain operations we do the following. 1) Consider only initial bases of families which are oriented open sets in Euclidean balls. 2) We choose a cocycle Poicar\'e dual to the diagonal in $M \times M$ supported in a small neighborhood of the diagonal. 3) We pull back this dual cocycle to the Euclidean balls by the evaluation maps just discussed. 4) When critical levels coalesce we form the Cartesian or tensor product of the dual cocycles and pull them back to the base of the input boundary family. 5) The pair consisting of the original base family and the pull back cocycle is the output chain. More formally the output chain is the the cap product of the input chain with the pulled back dual cocycle, $\sigma \cap U$ where $U$ is the pulled back cocycle. 6) If the input chain $\sigma$ is already a cap product pair $(\sigma' , U') \sim (\sigma' \cap U')$ then $(\sigma', U', U)$ will correspond to $(\sigma' \cap U') \cap U = \sigma' \cap (U \cap U')$. 

With this formulation using the Poicar\'e dual cocycle, we can enter the forbidden region beyond the open moduli space cut off by all essential circuits $\geq \epsilon$. One merely pulls back the products of dual cocycles via the specializing evaluation mappings. We have seen above in section 1.4 how single loop degenerations are dealt with and filled in. Multiple loop degenerations work very well because there are more equations than unknowns and the corresponding dual cocycle product has too large a dimension and becomes identically zero under restriction and pull back (see section 1.4).

There is a contractibility property of dual cocycles to any diagonal. Any two differ by a coboundary in the pair $(\textit{neighborhood}, \del \textit{neighborhood})$. Any two such coboundaries differ by a coboundary, etc. In the construction above these coboundaries are added in to the geometric homotopies to create chain homotopies between slightly different geodesic arc constructions. We do illustrative examples in the sections 3.3 to 3.8.

These are the elements of the proof of the main workhorse theorem. Let $(X,M)_*$ denote the singular chains in the equivariant mapping space $(X,M)$. See section 3.3 and 3.4 for the definition of $(X,M)_*$.

\begin{center}
\begin{figure}[h]
\centering
\includegraphics[scale=0.9]{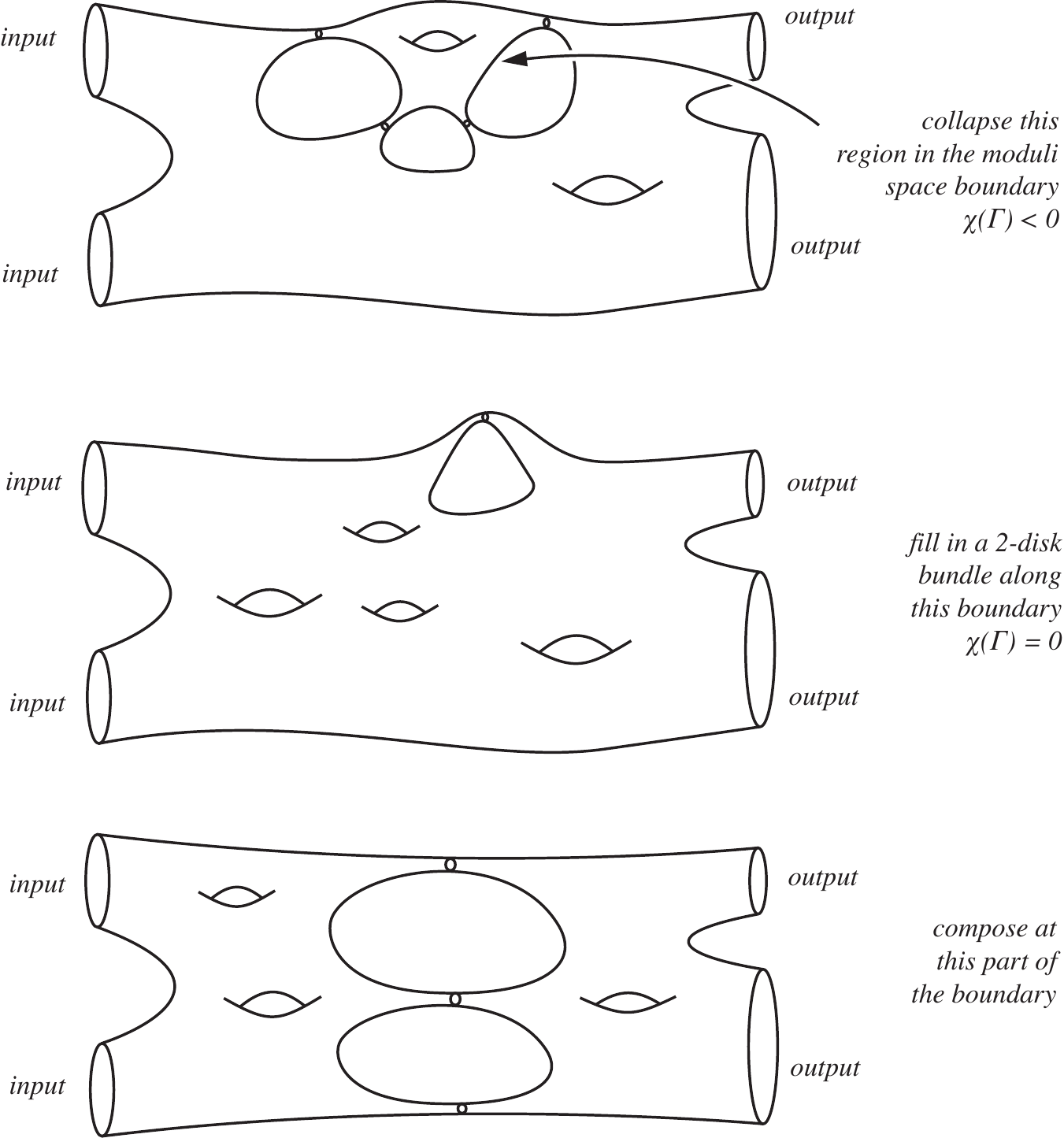}
\caption{Type A filling in plus composition boundary}
\label{fig:collapse-workhorsethm}
\end{figure}
\end{center}

\begin{center}
\begin{figure}[h]
\centering
\includegraphics[scale=1]{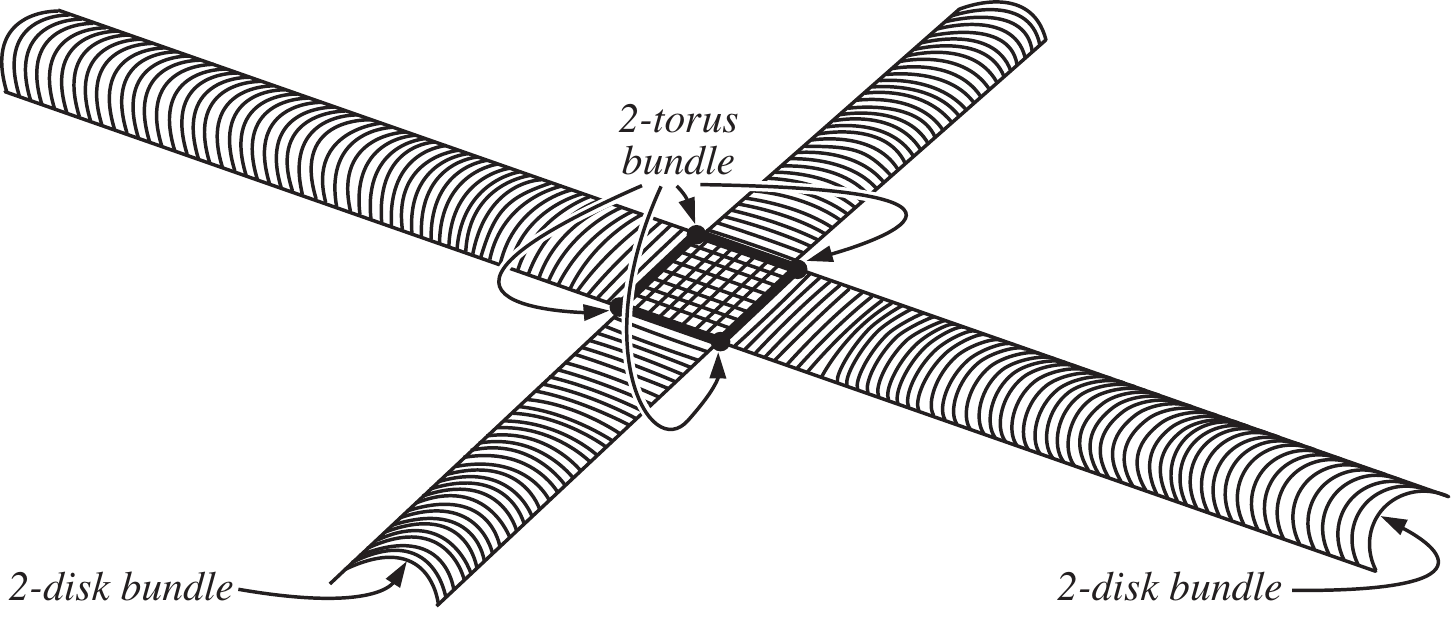}
\caption{}
\label{fig:2disk2torusbundle}
\end{figure}
\end{center}

\begin{thm}[Workhorse theorem]
For each $(k,l,g)$ there is a chain mapping of the equivariant chains on the compactified moduli space $\widehat{\mathscr{M}}(g,k,l)$ into the chain complex $Hom((I,M)_*,(\Sigma,M)_*)$. Here $\widehat{\mathscr{M}}(g,k,l)$ is the moduli space with part of its boundary, excluding the composition boundary, filled in or zipped up by coning off certain factors in fibered product decompositions of the non composition boundary. In more detail (see figure \ref{fig:collapse-workhorsethm}),

a) single circles that do not separate part of input from output or visa versa correspond to $S^1$ fibrations on the boundary. These are filled in with disk bundles.

b) collections of mutually nonhomotopic circles that do not separate correspond to torus bundles and get filled in by the intersection of cases of type a). See figure \ref{fig:2disk2torusbundle}

c) if a collection of mutually nonhomotopic circles separates off a component that has no input or output boundary, this region is coned off and the string topology construction is set equal to zero because of the negative Euler characteristic  argument in section 1.4. 

d) if a collection of circles separates off a component with a non-empty but balanced weight distribution they are treated as in $a)$ or $b)$ and filled in as circle bundles.

e) if a collection of circles separates off a component with an unbalanced set of weights, a composition is formed with the heavier part being the input or output of the composition depending on sign (see examples in section 3.7).

f) there is one more piece of noncomposition boundary -- the``outgoing lengths simplex boundary'' which is filled in in the closed string or equivariant case algebraically by working in the quotient by small loops in $M$.
\end{thm}

\begin{remark}
It seems likely the homology of these spaces $\widehat{\mathscr{M}}(g,k,l)$ can be known with present technology -- perhaps more easily than the open $\mathscr{M}(g,k,l)$. See the CUNY thesis of Fereydoun Nouri \cite{NouriThesis}.
\end{remark}

\subsection{The string bracket and jacobi relation at the chain level}

The Goldman bracket for curves up to homotopy on a surface generalizes to the string bracket operation on the equivariant chains $\lL^{S^1}_*(2) \stackrel{[,]}{\rightarrow} \lL^{S^1}_*(1)$ of the free $k$-loop space $k=1,2,3, \dots$ of an oriented $d$-manifold $M^d$. See definition 1 iii) for the definition of $\lL^{S^1}_*(k)$.

More exactly, there are chain mappings $\lL^{S^1}_*(2) 
\stackrel{I}{\leftarrow}
P(2,1)^G_* \stackrel{0}{\rightarrow} \lL^{S^1}_*(1)$ where $P(2,1)^G_*$ are the equivariant chains on the space of maps of a pair of pants $P(2,1)=\{S^2 - 3 \textrm{ disks} \}$ into $M^d$, $I$ and $O$ are restriction mappings to the input boundary (two components of the boundary of $P(2,1)$) and the remaining output boundary of $P(2,1)$. For a pair of pants $P(2,1)$ with $3$ boundary circles, two of which are input and one of which is output, the structure group $G$ is diffeomorphisms of $P(2,1)$ which are rotations on each boundary circle. The string topology construction produces maps $\lL^{S^1}_*(2) \stackrel{ST}{\rightarrow} P(2,1)^G_*$, and then $\lL^{S^1}_*(2) \stackrel{[,]}{\rightarrow} \lL^{S^1}_*(1)$ is the composition $O \circ ST$. $ST$ was discussed briefly above and will be discussed again momentarily.

\begin{defi}{($G$-equivariant chains on maps of a $G$-space $X$ into $M$)}
\begin{enumerate}
\item[i)]
First, consider a standard version of the equivariant singular chain complex. Here a $k-$simplex is a pair ($X$-bundle with specified structure group $G$ over the standard $k-$simplex, map of the total space into $M$) up to equivalence, where equivalence is an $X$-bundle isomorphism over the simplex satisfying: the obvious diagram of maps commutes.

\item[ii)]
Second, consider the diffuse equivariant version. Here we replace the standard $k$-simplex in i) by a pair $(s(\tilde{U}), \tilde{U})$ consisting of an open set $s(\tilde{U})$ in a $k+d-$simplex and a singular $d$-cochain $\tilde{U}$ whose closed support is in $s(\tilde{U})$. Then the bundle, the map and the equivalences need only be defined over the open set $s(\tilde{U})$.

The boundary map for i) is the usual one, and the boundary map for ii) is the direct analogue of the usual one: the alternating sum of the restrictions of the cocycle to the $k+d-1-$faces plus another term $\pm (s(\tilde{U}), \delta \tilde{U})$.

\item[iii)]
By $\lL^{S^1}_*(k)$ we mean the ``diffuse equivariant'' chains for the mapping space \{$k$-labeled circles, M\} with the structure group the $k$-torus acting by rotations on the $k$-labeled domain circles.

\item[iv)]
Later on we add a further equivalence relation allowing non identity diffeomorphisms on the base of the family.

\end{enumerate}
\end{defi}

Consider a bundle whose fiber is two labeled circles, with structure group the 2 torus, over $\Delta_k$, the standard $k$-simplex. Let $E(2)=$\{ordered pairs of points, one in each circle\}. Then given a map of the total space $E$ into $M$, form the map $E(2) \rightarrow M \times M$ by evaluating the map of $E$ into $M$. Pull back an apriori chosen cocycle $U$, Poincar\'e dual to the diagonal and supported on a small neighborhood of the diagonal, to get a cocycle $\tilde{U}$ on $E$. For each point in the support of $\tilde{U}$ the image point in $M \times M$ is near the diagonal by definition. Connect this pair of points in $M$ by a canonical short path (like the geodesic in some apriori chosen metric). It follows that for each point $p$ of the support of $\tilde{U}$ we get a map of a graph into $M$. The graph is made out of the two circles mapping into $M$, together with the short path associated to $p$

 \begin{figure}[h]
\centering
\includegraphics[scale=0.75]{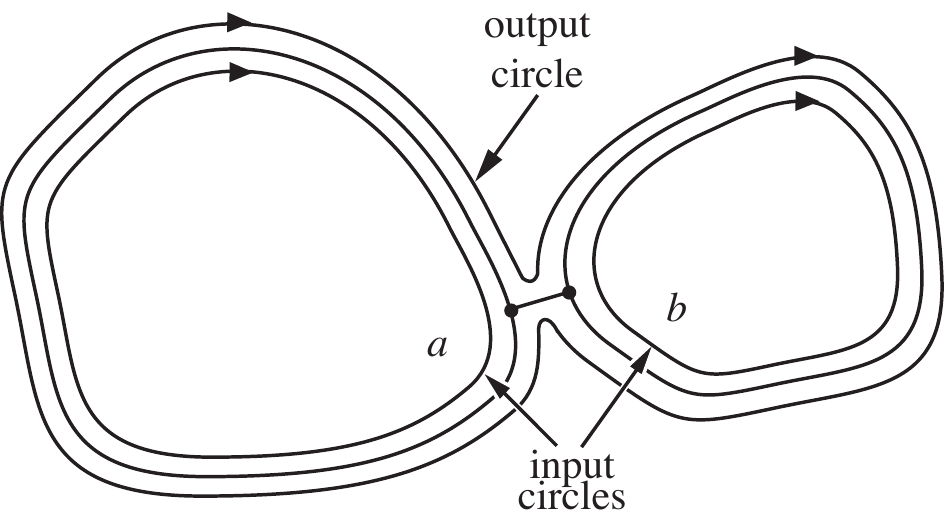}
\caption{The bracket}
\label{fig:bracket}
\end{figure}

The varying graph may be viewed as a varying spine of a varying $P(2,1)$ over $s(\tilde{U})$, the support of $\tilde{U}$. The pair $(s(\tilde{U}), \tilde{U})$ may be considered to be the output chain of $ST$. One can extend this construction from standard $k$-chains to diffuse equivariant chains using the cup product of cochains, namely $(s(V), V)$ yields $(s(V \cup \bar{U}), V \cup \bar{U})$ . 

The composition $O \circ ST : \lL^{S^1}_*(2) \rightarrow \lL^{S^1}_*(1)$ is the chain level bracket generalizing Goldman's.

\begin{center}
\begin{figure}[h]
\centering
\includegraphics[scale=0.75]{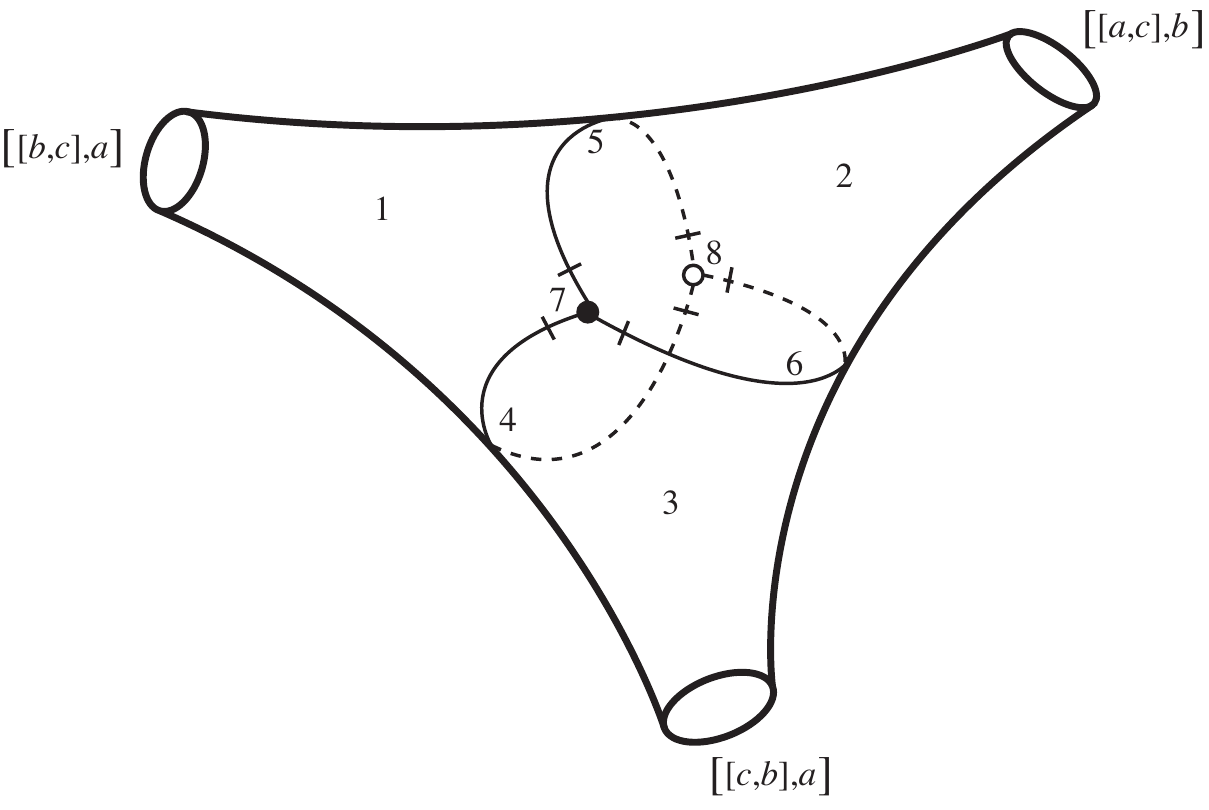}
\caption{Stratified cutoff moduli space of \{$S^2 - $4 points\}}
\label{fig:jacobimodspace}
\end{figure}
\end{center}

To discuss the analogue of the jacobi identity we consider $P(3,1)$, the two-sphere minus four disks with three labeled boundary components designated as input boundary and the remaining one as output. The moduli space of the two-sphere minus four points has three points at infinity. We will use this moduli space cut off near infinity to build out of 17 pieces a chain homotopy for the analogue of the jacobi relation for the string bracket. Each boundary component of the moduli space corresponds to a term in the jacobi identity.

\begin{figure}[h]
\centering
\includegraphics[scale=0.75]{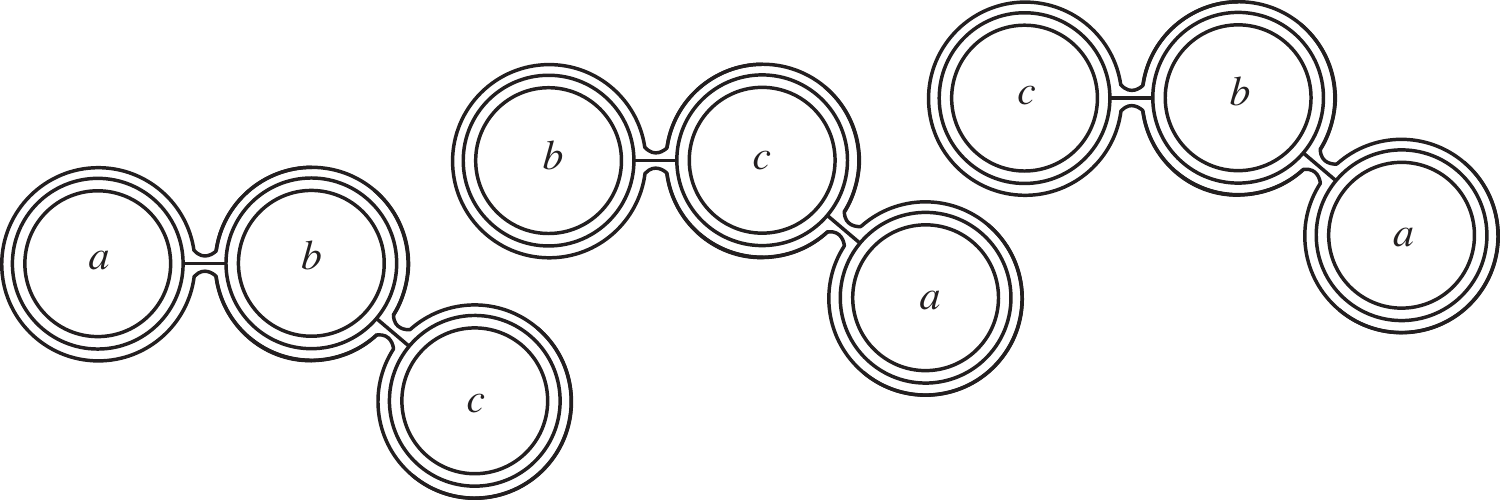}
\caption{}
\label{fig:ribbonrow}
\end{figure}

\begin{figure}[h]
\centering
\includegraphics[scale=0.75]{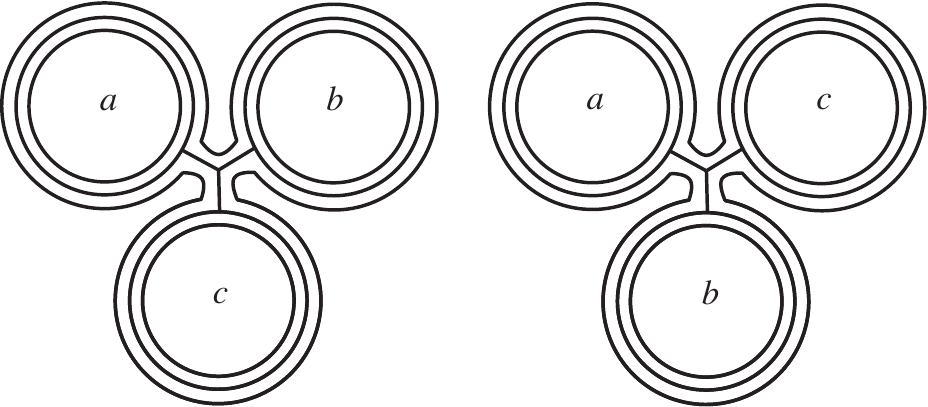}
\caption{}
\label{fig:ribbontripod}
\end{figure}

The moduli space is stratified as indicated in the figure \ref{fig:jacobimodspace}. The graphs of figures \ref{fig:ribbonrow} and \ref{fig:ribbontripod} determine ribbon graphs which thicken to Riemann surfaces. The input circles have equal length. Each stratum labeled will correspond to other chain operations $\lL^{S^1}_*(3) \to \lL^{S^1}_*(1)$ built according to string diagrams via diffusion intersection and short geodesic connections. Thus, strata 4, 5, 6 correspond to the three diagrams of figure \ref{fig:ribbonrow} depending on which circle is in the  middle. 

Strata 7, 8 correspond to the two diagrams of figure \ref{fig:ribbontripod} depending on the cyclic order of the labeled circles.

\begin{figure}[h]
\centering
\includegraphics{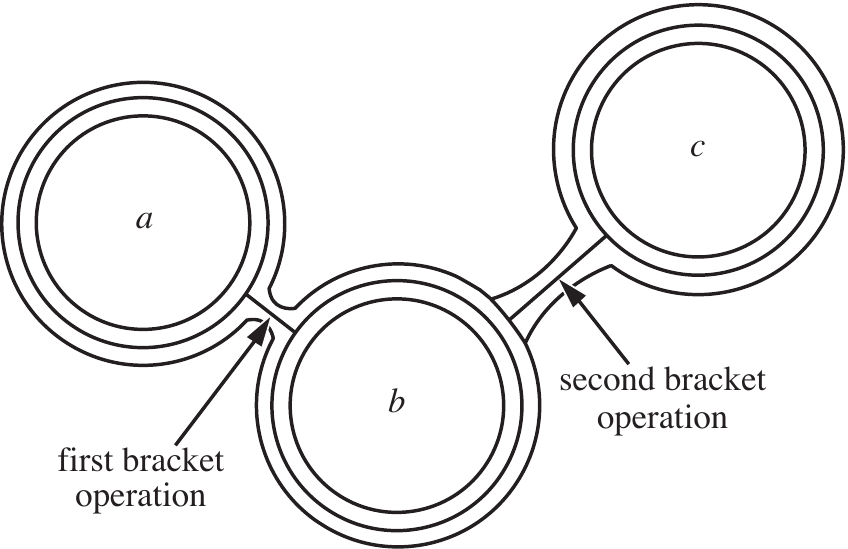}
\caption{$[[a,b],c]$}
\label{fig:ribboncomposition}
\end{figure}

Each of the boundary components being a composition is represented diagrammatically as in figure \ref{fig:ribboncomposition}, depending on labeling, e.g., which new circle enters into the second bracket operation.

The operations for strata 7 and 8 are defined by evaluating the simplex $\sigma$ in $\mathbb{L}^{S^1}_*(3)$ at three varying points, one on each component of the domain set of three circles to get a map 
\newline $T^3 \times \sigma \rightarrow M \times M \times M$. We pull back a Poincar\'e dual cocyle $U_{123}$ defined near the small diagonal $(x,x,x)$ to obtain the (base space, cocycle) piece of our equivariant chain in $P(3,1)_*$. Over each point in the support, we construct inside $M$ using short geodesics the (graph-spine) of the (sphere - four disks) corresponding to the surface labeled by the appropriate point in the moduli space of figure \ref{fig:jacobimodspace}. One way to do this is first take the geodesic convex hull of these nearby points, second take the barycenter of this convex hull, and third connect this barycenter to each of the three points by a short geodesic to construct the little corolla  or tripod connecting the three circles. This  tripod union of the three circles is the graph spine.

A similar construction is used for the graph spines of figure \ref{fig:ribbonrow}. Now we evaluate the map at one point on each of the outer circles and at an ordered pair of points on the middle circle. Label these 1, 2, 3, 4 reading left to right. We obtain a map of $T^4 \times \sigma \rightarrow M^4$. We form the Poincar\'e dual cocycles $U_{12}$ and $U_{34}$ near the diagonals of $(M \times M)_{12}$ and $(M \times M)_{34}$ and we pull back $U_{12} \cup U_{34}$ to obtain the diffuse base of our equivariant chain. We construct the graph spine of figure \ref{fig:ribbonrow} using short geodesics as before to obtain an equivariant chain in $P(3,1)_*$.

Now we come to the assembly step of these different pieces. Notice that as one of the circular arcs on the middle circle of figure \ref{fig:ribbonrow} shrinks to a vanishingly small length,  that point on stratum (4, 5 or 6) tends to the stratum (7 or 8) depending on the cyclic order. The two operations do not quite fit together. However, the geometric discrepancy is carried by the convex hull of the three nearby points. So it is easy to find a  homotopy of geodesics reconciling the slight difference. There is also a discrepancy between the diffusing classes used. We can argue abstractly as follows. The difference is carried by the small neighborhood of the diagonal and the difference is the coboundary of some $c$ there. Note again that the space of Poincar\'e dual cocycles is algebraically contractible in the sense that  two differ by a coboundary, and two such coboundaries differ by a coboundary, etc. We then combine the geometric homotopy with $c$ to obtain a chain homotopy between the two maps, reconciling the discrepancy. We put these homotopies over the small interval between the appropriate triple point and hash mark in figure \ref{fig:jacobimodspace}.)

We proceed to the last step. The composition constructions at the boundary of the cutoff moduli space of $\{S^2 - 4 \textrm{ points}\}$ corresponding to the geometric spine of figure \ref{fig:ribboncomposition} also differ slightly from the constructions of figure \ref{fig:ribbonrow} and figure \ref{fig:ribbontripod} (together with the small geometric homotopy reconciling them). Again the difference is carried by the convex hull of 3 nearby points in $M$ and can be reconciled by a geometric homotopy. Also again the Poincar\'e dual cocycles can be reconciled by a coboundary $c$. Together, the geometric homotopy and $c$ can be combined to yield a chain homotopy which we use over regions 1, 2, 3 respectively of figure \ref{fig:jacobimodspace}. 

\begin{remark}
\begin{enumerate}
\item[i)]
There are natural cartesian product maps, $\lL^{S^1}_*(k) \otimes \lL^{S^1}_*(l) \rightarrow \lL^{S^1}_*(k+l)$, inducing isomorphisms on homology

\item[ii)]
Under the permutation of component circles the generalized bracket map \\
$\lL^{S^1}_*(2) \rightarrow \lL^{S^1}_*(1)$ is skew-symmetric when $d$ is even because the fundamental class of the $T^2$ factor is reversed. It is symmetric when $d$ is odd because the Poincar\'e dual cocycle to the diagonal in $M^d \times M^d$ is also reversed under permutation of the factors.

\item[iii)]
The map in i) is graded symmetric for $k=l=1$. So we have proved the

\end{enumerate}
\end{remark}

\begin{prop}
The induced string bracket maps $\lL^{S^1}_*(1) \otimes \lL^{S^1}_*(1) \rightarrow \lL^{S^1}_*(1) $, has degree $-d+2$ and is graded skew-symmetric for $d$ even and graded symmetric for $d$ odd. 
\end{prop}

\begin{remark}
The common domain for all the pieces of the construction of the jacobi identity homotopy is $\lL^{S^1}_*(1) \otimes \lL^{S^1}_*(1) \otimes \lL^{S^1}_*(1) $. Namely, on regions 4, 5, 6, 7, 8 of figure \ref{fig:jacobimodspace} the domain is $\lL^{S^1}_*(3)$. On regions 1, 2, 3 of figure \ref{fig:jacobimodspace} it is $\lL^{S^1}_*(2) \otimes \lL^{S^1}_*(1)$ where the second tensor factor is, in turn, each labeled circle in the input boundary of $P(3,1)$. Thus the common domain is $\lL^{S^1}_*(1)^{\otimes 3}$. With this understanding we have the
\end{remark}

\begin{thm}
$\lL^{S^1}_*(1)$ has the structure of a graded differential Lie algebra of degree $-d+2$ in the sense that the bracket satisfies jacobi up to a \textbf{contractible set} of homotopies.
\begin{proof}
We add to the above discussion proof of the theorem the following remark: the construction of the jacobi homotopy $\lL^{S^1}_*(1)^{\otimes 3} \rightarrow P(3,1)^{G}_* $ was carried by small convex sets in $M^d$ and by the contractible sets of Poincar\'e dual classes to the diagonals. Thus the set of these homotopies forms a contractible object.
\end{proof}
\end{thm}

We will continue to work with these ideas to construct maps $\lL^{S^1}_*(1)^{\otimes k} \rightarrow P(k,1)^{G}_*$ using the moduli spaces of the 2-sphere $- (k+1)$ points in order to construct the hierarchy of higher homotopies of a infinity Lie structure on $\lL^{S^1}_*(1)$ ($\qQ$ coefficients). Using the homotopy theory of such structures (see Appendix), one then obtains a infinity Lie structure on $\hH^{S^1}_*$ the equivariant homology of the free loop space of $M$. Since Lie infinity structures on a complex $A$ can be reassembled as coderivations of square zero on the free cocommutative coalgebra $\Lambda^{c}(A)$ we will have shown:

\begin{thm}
The chain level string bracket construction, together with the moduli space chain homotopies, yields a coderivation differential of degree $-1$, $d=  d_2 + d_3 + \dots$ on the free graded algebra generated by the equivariant rational homology $\hH^{S^1}_*$ of the free loop space of $M^d$, shifted by $-d+3$. The differential $d$ is well defined up to isomorphism homotopic to the identity. \cite{Su}
\end{thm}

\begin{question}
We now know that $d_2$ or, equivalently, the bracket on equivariant homology $\hH^{S^1}_*$ is preserved by a homotopy equivalence between closed manifolds \cite{CKS}. It is conjectured that the entire package of String Topology up to equivalence is not a homotopy invariant of closed manifolds \cite{Su2}. One may already ask whether the structure of the higher terms of the differential on the above Lie infinity structure $(\Lambda \hH^{S^1}_*, d)$, which are not covered by the current theorems on homotopy invariance, e.g., \cite{CKS} constitutes a diffeomorphism invariant that is not a homotopy invariant.
\end{question}

\noindent
\textbf{Possible answer based on \cite{HopkinsLurie}.} The zipping up procedure for $(k,l,g)=(k,1,0)$ is not required here because the entire boundary is composition boundary. Thus our string topology construction in this case (which was in fact the part presented homologically in \cite{CS}) is perhaps homotopy equivalent to a construction of Hopkins-Lurie. Their construction is a homotopy invariant \cite{HopkinsLurie}.

\subsection{The string cobracket at the chain level}

\begin{figure}[h]
\centering
\includegraphics[scale=0.8]{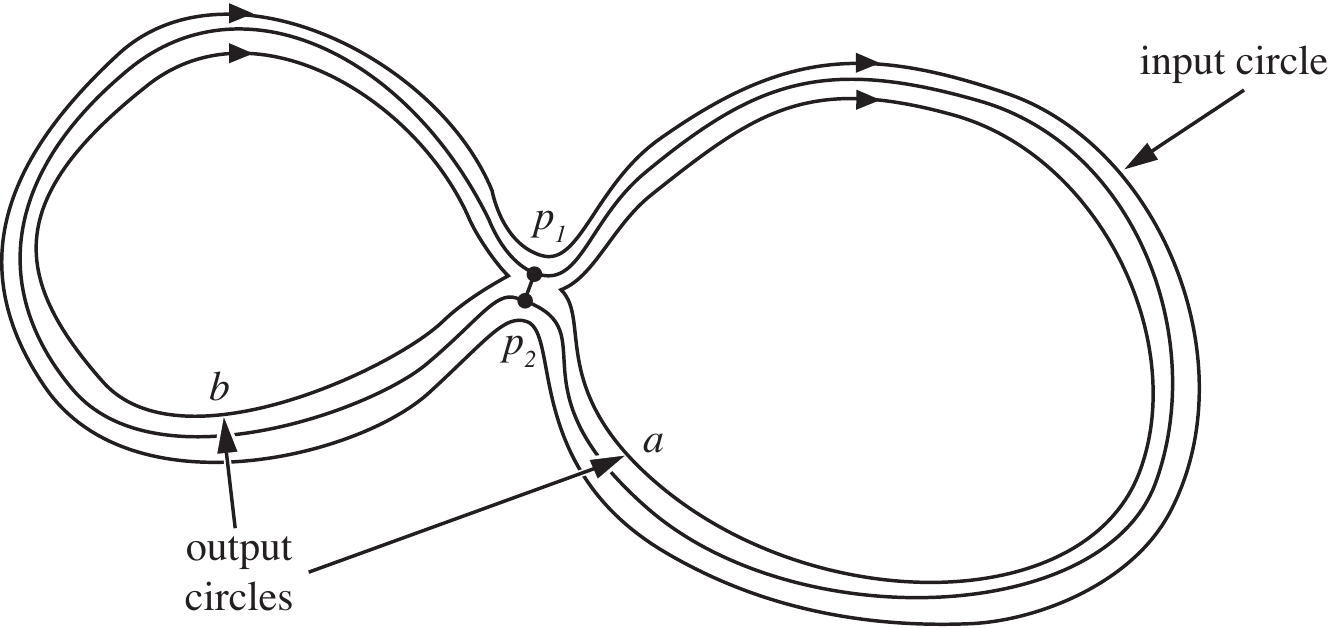}
\caption{The cobracket}
\label{fig:cobracket}
\end{figure}

Choose a chain in $\lL^{S^1}_*(1)$ (Definition 1 iii)) with base $B$ and total space circle bundle $E$. Over each point of $B$ put the configuration space $E(2)$ of ordered pairs of unequal points on the circle as a new fiber over this point of $B$. Compactify this fiber by blowing up the diagonal in the ordered pairs on the circle. Map this new total space into $M \times M$ by evaluating the given map of $E$ into $M$ at the various point pairs. We pull back the Poincar\'{e} dual cocycle defined in a small neighborhood of the diagonal. Over each point of the support of this pull back we have an ordered pair of nearby points in M. Connect these by a short geodesic and define a map of the spine shown in figure \ref{fig:cobracket} into M.

This defines an equivariant chain $T$ of maps of $P(1,2)$ into $M$. By restricting to the outer boundary of $P(1,2)$ we obtain two output circles. We order (and orient) the two output circles using the ordering of the point pair combined with the orientation of the input circle. For example, one way to order the circles is to say the first ouput circle is the one first traversed starting from the first point of the point pair and going in the direction of the orientation of the input circle. Note that blowing up the diagonal of the torus was \textbf{necessary} so that this ordering construction extends continuously to the compactification.
 
We assume by averaging the Poincar\'{e} dual cocycle is symmetric (up to sign, if $d$ is odd) under the flip of factors in $M \times M$. 
 
There is the involution of our chain corresponding to interchanging the order of the two points on the input circle. To understand this symmetry first assume the dimension, $d$, of $M$ is odd. Under the flip of factors in $M^d \times M^d$ the Poincar\'e dual cocyle is by assumption converted to minus itself. Also the orientation of the base of our chain is reversed. Thus for the diffuse chain pair $(s(\tilde{U}), \tilde{U})$ the flip of factors induces an orientation preserving symmetry of our chain.

Now we introduce a further equivalence relation on chains (foreshadowed in Definition 1 iv) and which we call attention to by using a prime superscript) that identifies two chains that are related by an orientation preserving diffeomorphism of the base of the chain. Thus our chain becomes two copies of the same reduced chain by choosing a fundamental domain of the involution. Moreover, this reduced chain is symmetric under the involution of $\lL^{S^1}_*(2)'$ interchanging the components for $d$ odd.
 
For $d$ even the reduced chain is reversed under the flip of components in $\lL^{S^1}_*(2)'$. There is a subtle sign issue to choose one representative of the reduced chain in the even case. 
Let us then study the boundary of this cobracket construction defined by the reduced chain. There are three possible contributions to the boundary coming from

 \begin{enumerate}
 \item[a)]
 the boundary of the two-point configuration space of the input circle.
 \item[b)]
 the boundary of a fundamental domain for the involution exchanging the order of the points in the two-point configuration. This fundamental domain is chosen to form representations of the reduced chain above.
 \item[c)]
 the usual boundary of the input chain for the cobracket construction.
\end{enumerate}

We will kill a) by working modulo small loops (see below). We have in effect killed b) by identifying to zero any chain which admits an orientation-reversing automorphism. This is the case for the new boundary of the fundamental domain of the involution, the base of the reduced chain. For c) we naturally do nothing.\\

Now we return to a). We work on the quotient of the equivariant chain complex by the ``$\epsilon$-small loops.'' Here $\epsilon$-small is defined by any number $\epsilon > 0$ so that for any two points of distance $\leq \epsilon$ in the apriori chosen metric there is a unique geodesic between them.\\

\begin{defi}
A chain in $\lL^{S^1}_*(k)'$ belongs to ($\epsilon$-small loop)$_*$ if there is a covering of the base of the chain so that in each open set of the covering there is an index $i$ so that the $i^{th}$ component circle has length $< \epsilon$ in the target.
\end{defi}

\begin{remark}
We will see below the bracket is still well defined on this quotient by $(\epsilon-small \; loops)_* $.
\end{remark}

\subsection{Cojacobi chain homotopy at the chain level}

\begin{figure}[h]
\centering
\includegraphics{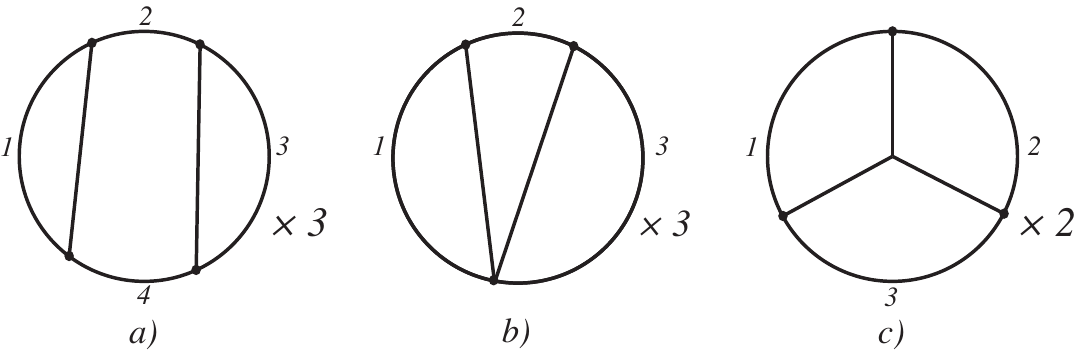}
\caption{8 string diagrams for the cojacobi chain homotopy}
\label{fig:cobracket^2}
\end{figure}

Now we discuss the analogous cojacobi chain homotopy for the cobracket. We use the moduli space for ($S^2- $4 points) again as a key cartesian factor for the combinatorial moduli space of $P(1,3)$, the combinatorial surface with one input circle and three output circles. In addition to the moduli space of figure \ref{fig:jacobimodspace} we have a circle factor for the marking on the input circle and the 2-simplex of parameters describing the distribution of combinatorial length on the three output circles. Thus ignoring the circle, our moduli space pseudomanifold with boundary pair is essentially the double suspension of the 2-disk minus two smaller disks. We see a 4-dimensional chain with three cycles on the boundary of dimension 3. The pieces organizing the cojacobi homotopy correspond to the eight string diagrams in figure \ref{fig:cobracket^2}. The same type of argument detailed above for jacobi will produce a chain mapping $\lL^{S^1}_*(1)'' \rightarrow P(1,3)^{G \prime \prime}_* $ so that the restrictions to the three cycles correspond to the three composition terms appearing in cojacobi. The double prime refers to modding out by small output loops and degenerate chains (see next section) as well as by the equivalence relation involving diffeomorphisms of the bases of chains of maps, which assures us that near the boundary of the 2-simplex of output length distribution we have the zero mapping. The extension to the interior provides the chain homotopy of the cojacobi expression to zero.

\subsection{Extending the bracket to the quotient by small strings}

The idea is that we could have worked from the beginning with $\lL^{S^1}_*(k)'$ modulo geometrically degenerate chains -- those which, as maps into some space, have a lower dimensional image than their domain dimension.

A chain obtained by bracketing with a family of constant closed strings has this property from the definition. This is because we bracket with every point of the circle mapping to the constant string. Now the subcomplex of small strings is chain equivalent to the subcomplex of constant strings. Combining these two facts gives the desired extension.

\subsection{The chain homotopy for drinfeld compatibility of the bracket and cobracket}

This homotopy arises from the moduli space of 2 inputs and 2 outputs for ($S^2-4$ points). Now besides the 2 dimensions of the  ($S^2-3$ disks) of figure \ref{fig:jacobimodspace}, we have the torus of marks on the input boundary and the interval of length distributions on the output boundary. As usual we ignore the 2-torus factor and concentrate on the rest, which is ($S^2-3$ disks) in effect suspended once by crossing with the interval and working mod the endpoints of the interval (corresponding to working in the quotient by small loops).

 \begin{figure}[h]
\centering
\includegraphics{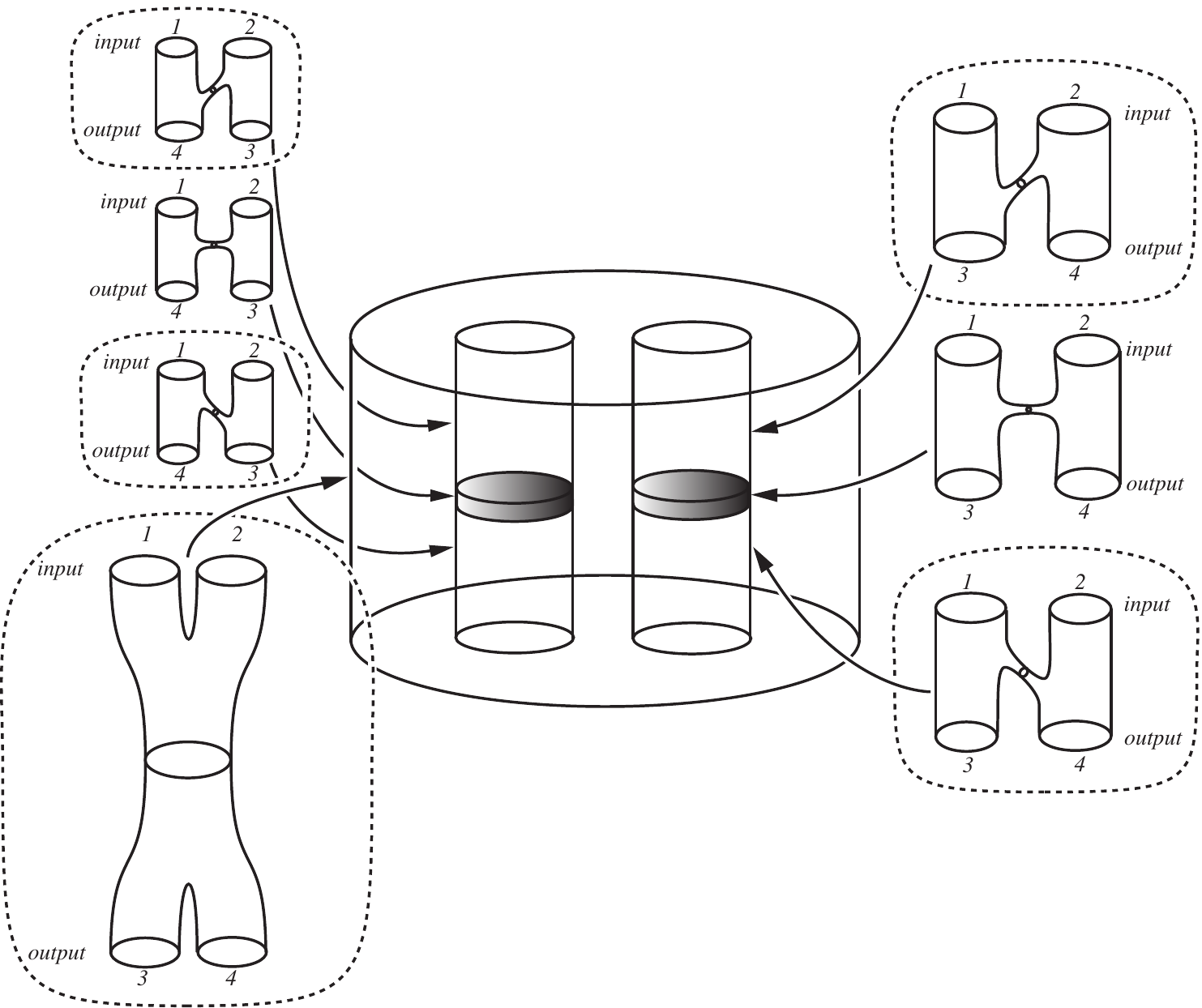}
\caption{Drinfeld compatibility five term relation (circled)}
\label{fig:drinfeld}
\end{figure}

In the figure \ref{fig:drinfeld} the outer boundary component corresponds to the composition: do the bracket then the cobracket. Each inner component corresponds to two terms in the composition boundary. The upper piece on the right side of the figure 
is the composition: cobracket of the input 2 followed by bracketing in the input 1 on the right factor. This occurs on the part of the interval where the length of the output 3 is greater than the length of the output 4. The lower piece on the right hand side 
corresponds to the composition: cobracket of the input 1 followed by the bracket of the input 2 on the left. Here output length 4 is greater than output length 3.

When the outputs' lengths are approximately equal the string topology construction yields a small image for the separating circle. This can be filled in by the discussion of section 1.4, $\chi(\Gamma) = 0,$ case. This filling in is indicated by the 2-handle in the figure.

The other inner boundary is treated similarly. In conclusion the construction extends over the moduli space with these two 2-handles added. Now the boundary has five components and each corresponds to a composition. This is the five-term relation of drinfeld compatibility.
$$\Delta[x,y] = [\Delta x, y] \pm [x, \Delta y] $$ 
\noindent
where each term on the right is two terms, the two terms in the definition of the action of a Lie algebra on its tensor square.

\subsection{The chain homotopy for the involutive identity}

Now we use our first higher genus moduli space, the space of all \{torus -- 2 points \}. The moduli space of the torus with one tangent marked point is the complement of the trefoil knot in $S^3$ -- by an argument usually attributed to Quillen. Adding in another puncture without a mark for the output is the 5-dimensional universal fibration with fiber ($T^2-$point) over the complement of the trefoil.

\begin{figure}[h]
\centering
\includegraphics{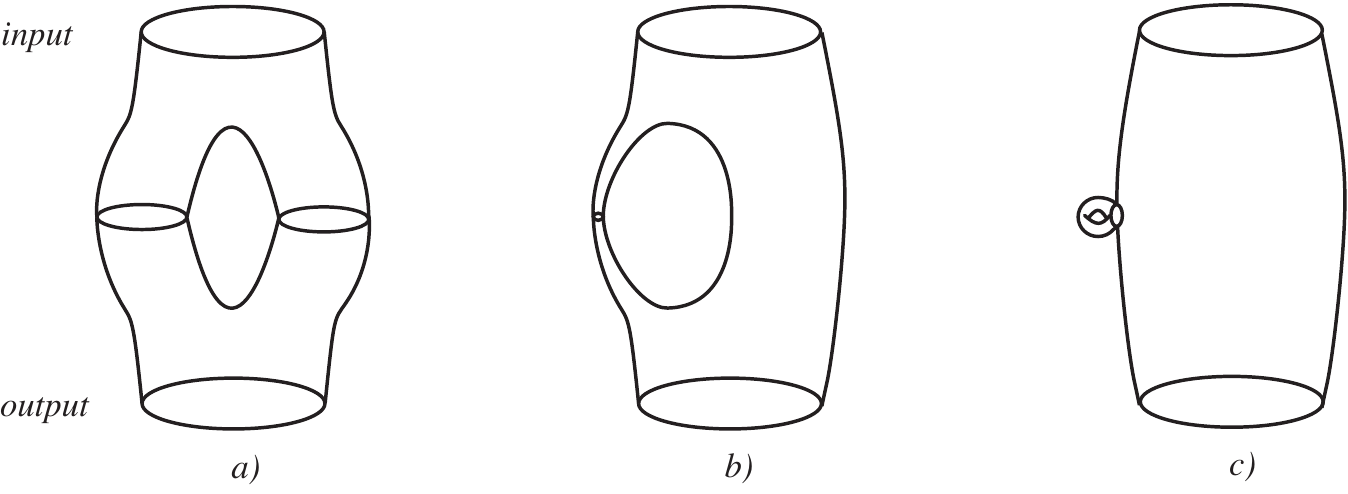}
\caption{Involutive relation - boundary terms}
\label{fig:involutivedegenerate}
\end{figure}

The composition boundary of this moduli space corresponds to two disjoint embedded curves which separate input from output. See figure \ref{fig:involutivedegenerate} a).

We will fill in the rest of the boundary. If only one curve is collapsing, there are two cases.

\begin{enumerate}
\item[i)]
The curve is essential in the torus and we fill in by the $\chi(\Gamma) = 0$ argument above. See figure \ref{fig:involutivedegenerate} b).

\item[ii)]
The curve separates the input and output from the rest of the surface and we fill in by the $\chi(\Gamma) < 0$ which implies the locus of the string topology construction is vacuous near this part of the boundary of moduli space. See figure \ref{fig:involutivedegenerate} c).

\end{enumerate}

If an additional curve is collapsing and we are not in the composition boundary we are contained in ii).

Thus after filling in we get a chain operation whose boundary is the composition of cobracket followed by bracket. This has been dubbed the involutive relation \cite{CS2}.

\appendix

\section*{Appendix: Homotopy theory of the master equation}

For simplicity, we restrict to linear and quadratic terms but there is no obstruction to treating the general case.

\begin{defi}[The master equation] Symbolically $dX + X \ast X  + LX = 0$ where $X = \{X_\alpha\}$ is a linear basis indexed by some indexing set \{$\alpha$\}, $X \ast X$ is a sum over 
a collection $\ast$  of binary operations combining the various basis elements $X_{\alpha}$, and $LX$ is a sum over a  collection of linear operations $L$ on the span of $\{X_{\alpha} \}$. This may be formalized one way in the language of universal algebras \cite{Bergman, Cohn}.
\end{defi}

Actually, we think of such a master equation as a presentation of a free differential algebra generated by the $X_{\alpha}$ with the differential on the generators defined by the master equation.
The condition $d \circ d  = 0$ should be a formal identity in the free algebra generated by the $\{X_{\alpha} \}$. The algebra is meant to be free as an algebra over an \textit{operad}  $O$ generated by the binary operations in $\ast$ and the unary operations of $L$. These operations may satisfy relations like jacobi  so that $O$ itself need not be free. We further assume for what follows that there is a partial ordering of the indexing set, with all descending chains finite, so that the right hand side of the equation for $dX_\alpha$ only contains variables of strictly lower index.  This property of a free dga is called the \textit{triangular property}.

Here are some interesting examples of non free differential graded algebras over relevant gluing or composition $O$'s. Let $C$ denote a chain complex.
\begin{enumerate}
\item
$Hom'(C) =  \oplus_{k > 0} Hom (C^{\otimes k} , C)$ where binary $\ast$ operations are obtained by substituting the output  for some $k$ into one of the inputs for some $j$. There are $j$ such binary operations given $k$ and $j$. Note the Leibniz rule holds for the binary operations relative to the natural differentials on the hom spaces.
\item
$Hom''(C) =  \oplus_{k > 0, j > 0} Hom (C^{\otimes k} , C^{\otimes j})$  where binary $\ast$ operations are substituting an output for some $k, j$ into one of the inputs for some $k', j'$ . There are in this case $jk'$ binary operations. Again Leibniz holds.
\item
$Hom'''(C) = \oplus_{k > 0, j > 0} Hom(C^{\otimes k} , C^{\otimes j})$ where binary $\ast$ operations are substituting one or more outputs for some $k, j$ into a set of inputs for some $k',j'$ . 
Leibniz holds here as well.
\item
$Hom ^{(IV)} =$ same spaces and operations as 3) with tensor products thrown in. Leibniz holds here.
\end{enumerate}

These examples, while not free, admit maps from free triangular dgas  (in each context, i.e., for each dga over an operad $O$) inducing isomorphisms on homology. This is true for arbitrary dgOas i.e., differential graded algebras over the operad $O$. Such dgOa maps are called \textit{resolutions}. There is a notion of homotopy for maps from free triangular dgOas into arbitrary dgOas allowing the definition of homotopy equivalence of free triangular dgOas.
Resolutions are well defined up to homotopy equivalence and the  homotopy equivalence is well defined up to homotopy. The latter uses a lifting proposition that says maps from free triangular dgOas into arbitrary dgOas can be lifted through dgOa homology isomorphisms.
$$
\xymatrix{
&
&
R \ar[d]^{\simeq_H} \\
&
F \ar[r]^S \ar@{-->}[ru]
&
A
}
$$

If $C$ and $D$ are chain equivalent chain complexes  one can also show, for each of the cases 1), 2), 3) the respective resolutions of $Hom(C)$ and  $Hom(D)$ are homotopy equivalent.

\begin{defi}
A ``master equation package'' is a triple $(F, S, A)$ where $F$ and $A$ are dgOa algebras with $F$ free and triangular and $S: F \rightarrow A$ a dgOa map. Two master equation packages $(F, S, A)$ and $(F', S', A')$ are homotopy equivalent if there are homotopy equivalences $f$ and $r$ between $F$ and $F'$ and $R$ and $R'$ respectively so that the obvious diagram, below, using liftings is homotopy commutative. Here $R$ and $R'$ denote triangular free resolutions of $A$ and $A'$ respectively.
$$
\xymatrix{
F\ar[r]^S \ar[d]^f
&
A
&
R\ar[l]_{\simeq_H} \ar[d]^r \\
F' \ar[r]^{S'}
&
A'
&
R'\ar[l]_{\simeq_H}
}
$$

\end{defi}

\begin{enumerate}
\item[(a)]
If $A$ is one of the examples 1), 2), 3), 4) we say the master equation package defines an \textit{infinity algebraic structure} on $A$ of the \textit{form given by the homotopy type} of $F$.
\item[(b)]
If $S$ is a homology isomorphism then we have called the master equation package a resolution. In this case the dgOa  is $A$ itself.
\item[(c)]
If $F$ is derived from looking at a system of moduli spaces and their codimension one frontiers, $A$ is a chain complex of geometric objects where the moduli spaces live, $A$ is provided with operations needed to describe the codimension one frontiers , and S is defined using the moduli spaces we say the master equation package arises from the compactness and gluing theory of a system of moduli spaces.
 \end{enumerate}

In the text, the string topology for closed strings is a package of type (a) (3) where $C$ is the homology of equivariant chains on the loop space of a manifold modulo constant loops. We have called the algebraic structure arising there a quantum lie bialgebra which incorporates the involutive identity up to a first homotopy at least, (section 3.6) and the rest of the the three lie bialgebra identities up to infinite homotopy.
The open string topology was presented as an example of type (c) where $A$ had the form of a hom of  chains on maps of strings into M into  equivariant chains of maps of surfaces into M.

For applications to symplectic topology, $F$ will be the recipe arising in the compactness and gluing pictures, $S$ will be defined by the moduli spaces to first approximation and $A$ will be defined by thinking of families of solutions to the elliptic equations as chains in algebraic topology and providing $A$ with operations for the gluing. 

If some of these operations involve transversality in the first approximation, i.e., they are like string topology operations, then higher order approximations to success will involve using the device of type (a) packages to really define the operations on $A$ completely.

For example, if a string bracket is needed for a term in the frontier of the moduli space one builds a lie infinity structure as in string topology using the top cells of a genus zero moduli space construction. Then one, to a higher order approximation, defines $S$  using  this lie infinity structure instead of  just the transversally defined string bracket. There are complicated details here but I suppose this perspective may also help out with internal transversality issues in the elliptic equations themselves.
For example there is a nice language of Kuranishi and emphasized by Fukaya, Oh, Ohta and Ono for describing a locally finite  infinite system of finite dimensional problems \cite{FOOO} to which the infinite dimensional elliptic problem reduces. One can imagine applying the Poincar\'e dual cocycle approach to this infinite system of finite dimensional problems, and then building in the homotopies to heal discrepancies as was done above for string topology.


\begin{thebibliography}{9}

\bibitem{Abbaspour} Abbaspour, H., \textit{On String Topology of 3-manifolds}, PhD dissertation, City Univ. of New York, 2003.
%
\bibitem{ACG} Abbaspour, H., Cohen, R.L., Gruher, K., \textit{String topology
of Poincar\'e duality groups}, arXiv:math/0511181.
%
\bibitem{AbbondandoloSchwarz}
Abbondandolo, A. and Schwarz, M., \textit{
On the Floer homology of cotangent bundles},
Comm. Pure Appl. Math. 59 (2006), no. 2, 254-316.
%
\bibitem{AndersenUeno}
Andersen J. E. and K. Ueno, K.,  \textit{ Geometric construction of modular functors from conformal field theory}, Journal of Knot theory and its Ramifications 16  (2007), no. 2, 127 Ð 202.
%
\bibitem{BV1}
Batalin, I., Vilkovisky, G., \textit{Gauge algebra and quantization}, Phys. Lett. 102B, 27 (1981).
%
\bibitem{BV2}
Batalin, I., Vilkovisky, G. \textit{Quantization of gauge theories with linearly dependent generators},
Phys. Rev. D29, 2567 (1983).
%
\bibitem{Bergman}
Bergman, G., \textit{An Invitation to General Algebra and Universal Constructions}, Henry Helson, Berkeley, 1998.
%
\bibitem{Bodigheimer} B\"odigheimer, C-F., \textit{Configuration Models for Moduli Spaces of Riemann Surfaces with Boundary}, Abh. Math. Sem. Univ. Hamburg (2006), no. 76,191-233.
%
\bibitem{Chas}
Chas, M., {\it Combinatorial Lie bialgebras of curves on surfaces}, 
Topology 43 (2004), no. 3, 543-568.
%
\bibitem{Chas2}
Chas, M., \textit{Minimal intersection of curves on surfaces}, arXiv:0706.2439.
%
\bibitem{CS}
Chas, M. and Sullivan, D., {\it String Topology}, arXiv: math.GT/9911159.
%
\bibitem{revisedCS}
Chas, M. and Sullivan, D., {\it Compactified String Topology}, in progress.
%
\bibitem{CS2} Chas, M. and Sullivan, D., \textit{Closed string operators in
topology leading to Lie bialgebras and higher string algebra},  The legacy
of Niels Henrik Abel,  771--784, Springer, Berlin, 2004.
%
\bibitem{CK}
Chas, M., and Krongold, F.,  {\it An algebraic characterization of simple closed curves}, 
preprint.
%
\bibitem{Chataur} Chataur, D., \textit{A bordism approach to string
topology}, Int. Math. Res. Not. (2005), no. 46, 2829-2875.
%
\bibitem{ChenThesis} 
Chen, X., \textit{On a General Chain Model of the Free Loop Space  
And String Topology}, PhD dissertation, State Univ. of New York, Stony Brook, 2007.
%
\bibitem{CieliebakLatschev}
Cieliebak, K. and Latschev, J., \textit{The role of string topology in symplectic field theory}, arXiv:0706.3284.
%
\bibitem{FRCohen}
Cohen, F. R., \textit{
Cohomology of braid spaces},   Bull. Amer. Math. Soc.  79(1973), 763-766.
%
\bibitem{Cohen}  Cohen, R. L.,  \textit{Morse theory, graphs, and string
topology},  Morse theoretic methods in nonlinear analysis and in symplectic
topology,  149--184, NATO Sci. Ser. II Math. Phys. Chem., 217, Springer,
Dordrecht, 2006.
%
\bibitem{CG}
Cohen, R. L. and Godin, V., {\it A polarized view of string topology} in {\it Topology, Geometry and Quantum Field Theory: Proceedings of the 2002 Oxford Symposium in Honour of the 60th Birthday of Graeme Segal} London Mathematical Society Lecture Note Series, Cambridge: Cambridge University Press, 2004.
%
\bibitem{CJ} Cohen, R. L. and Jones, J. D. S., \textit{A homotopy theoretic
realization of string topology},  Math. Ann.  324  (2002),  no. 4, 773--798
%
\bibitem{CohenJonesYan}
Cohen, R. L., Jones, J. D. S. and Yan,  J., \textit{The loop homology algebra of spheres and projective spaces}, Progr. Math., vol. 215, Birkhauser, Basel, 2003, 77-92.
%
\bibitem{CKS}
Cohen, R. L., Klein, J. and Sullivan, D., {\it The homotopy invariance of the string topology loop product and string bracket}, preprint, submitted to Journal of Topology, arXiv:math/0509667v1.
%
\bibitem{CohenVoronov}
Cohen, R. L. and Voronov, A., \textit{Notes on string topology}, String topology and cyclic homology, 1--95, Adv. Courses Math. CRM Barcelona, BirkhŠuser, Basel, 2006 and arXiv:math.GT/0503625.
%
\bibitem{Cohn}
Cohn, P. \textit{Universal Algebra}, Springer, Berlin, 1981.
%
\bibitem{Connes}
Connes, Alain,  \textit{Noncommutative differential geometry}, Inst. Hautes ƒtudes Sci. Publ. Math. (1985), no. 62, 257-360.
%
\bibitem{Costello}
Costello, K., \textit{Renormalisation and the Batalin-Vilkovisky formalism}, 2007 arXiv:0706.1533.
%
\bibitem{Costello2}
Costello, K., \textit{Topological conformal field theories and Calabi-Yau categories}, 2004.
arXiv:math/0412149. Advances in Mathematics,  Volume 210, Issue 1, 
March 2007.
%
\bibitem{DeuschelStroock}
Deuschel, J-D. and Stroock, D., \textit{ Large deviations}, 
Pure and Applied Mathematics, 137. Academic Press, Inc., Boston, MA, 1989. 
%
\bibitem{EkholmEtnyreNgSullivan}
Ekholm, T., Etnyre, J., Ng., L. and Sullivan, M., \textit{in progress}.
%
\bibitem{EliashbergGiventalHofer}
Y. Eliashberg, A. Givental and H. Hofer, 
\textit{Introduction to symplectic field 
theory}, Geom. Funct. Anal. Visions in Mathematics special volume, part II, (2000), 560Ð 
673. 
%
\bibitem{FOOO}
Fukaya, K., Oh, Y-G., Ohta, H. and Ono K., \textit{Lagrangian intersection Floer theory--anomaly and obstruction}, preprint, http://www.kusm.kyoto-u.ac.jp/fukaya/fukaya.html, 2000. 
%
\bibitem{GadgilNg}
Gadgil, S., and Ng. L., \textit{The Chord algebra and fundamental groups}, appendix to \textit{Knot and braid invariants from contact homology II} by Lenny Ng, Geom. Topol. 4 (2005), 1603-1637.
%
\bibitem{Gerst}
Gerstenhaber, M., {\it The cohomology structure of an associative ring},
The Annals of Mathematics, 2nd Ser., Vol. 78, (1963), no. 2, 267-288.
%
\bibitem{Getzler}
Getzler, E., \textit{Batalin-Vilkovisky algebras and two-dimensional topological field theories}, Comm. Math. Phys. 159 (1994), no. 2, 265-285.
%
\bibitem{GetzlerJonesPetrack}
Getzler, E., Jones, J., Petrack, S., \textit{Loop spaces, cyclic homology and the Chern character}, Operator algebras and applications, Vol. 1, 95--107, 
London Math. Soc. Lecture Note Ser., 135, Cambridge Univ. Press, Cambridge, 1988. 
%
\bibitem{Godin} Godin, V., \textit{Higher string topology operations}, to appear.
%
\bibitem{Gol} Goldman, W.,
{\it Representations of fundamental groups of surfaces}, 
Springer Lecture Notes in Math., 1167 (1985).
%
\bibitem{Gol2} Goldman, W.,
{\it Invariant functions on Lie groups and
Hamiltonian flows of surface group
representations}, Invent. Math. 85 (1986),
no. 2, 263-302.
%
\bibitem{GS} Gruher, K. and Salvatore, P., \textit{Generalized string topology
operations}, arXiv:math/0602210.
%
\bibitem{Hochschild}
Hochschild, G., \textit{On the cohomology theory for associative algebras}, Ann. of Math. (2) 47, (1946). 568-579.
%
\bibitem{HopkinsLurie} Hopkins, M. and Lurie, J., \textit{Topological Quantum Field Theories in Low Dimensions}, in preparation.
%
\bibitem{Hu} Hu, P., \textit{Higher string topology on general spaces},
Proc. London Math. Soc. (3)  93  (2006),  no. 2, 515-544.
%
\bibitem{Jaco}
Jaco, W., {\it Heegaard splittings and splitting homomorphisms}, Trans. Amer. Math. Soc. 144 (1969), 
365Ð379.
%
\bibitem{Jones} Jones, J. D. S., \textit{Cyclic homology and equivariant homology}, Invent. Math. 87 (1987), no. 2, 403Ð423.
%
\bibitem{Kaufmann} Kaufmann, R., \textit{A proof of a cyclic version of
Deligne's conjecture via Cacti}, arXiv:math/0403340.
%
\bibitem{Kaufmann2}
Kaufmann, R., \textit{Moduli space actions on the Hochschild co-chains of a Frobenius algebra I: cell operads}, Journal of Noncommutative Geometry 1 (2007), no. 3, 333--384, arXiv:math/0606064.
%
\bibitem{Kaufmann3}
Kaufmann, R.,\textit{Moduli space actions on the Hochschild co-chains of a Frobenius algebra II: correlators}, arXiv:math/0606065.
%
\bibitem{KP} Kaufmann, R. and Penner, R. C., \textit{Closed/open string
diagrammatics},  Nuclear Phys. B  748  (2006),  no. 3, 335-379.
%
\bibitem{Kontsevich} Kontsevich, M., \textit{Deformation quantization of Poisson Manifolds, I.} Lett.Math.Phys. 66 (2003) 157-216., arXiv:q-alg/9709040. 
%
\bibitem{Kontsevich2}
Kontsevich, M. and Soibelman, Y., \textit{Notes on A-infinity algebras, A-
infinity categories and non-commutative geometry}, {I}, (2006) 
arXiv:math/0606241.
%
\bibitem{Kontsevich3}
Kontsevich, M., Lecture at the Hodge centennial conference, Edinburgh, 
2003,  Lectures in Miami, 2004.
%
\bibitem{Loday}
Loday, Jean-Louis, \textit{Cyclic Homology}, Grundlehren der mathematischen Wissenschaften Vol. 301, Springer (1998).
%
\bibitem{Mandell1}
Mandell, Michael M., \textit{Equivariant $p$-adic homotopy theory}, Topology Appl. 122 (2002), no. 3, 637-651.
%
\bibitem{Mandell}
Mandell, Michael M., \textit{$E_\infty$ algebras and $p$-adic homotopy theory}, Topology 40 (2001), no. 1, 43-94.
%
\bibitem{McCrory}
McCrory, Clint G., \textit{Poincar\'e Duality in Spaces with Singularities}, Thesis, Brandeis, 1972.
%
\bibitem{Menichi} Menichi, L., \textit{Batalin-Vilkovisky algebras and cyclic
cohomology of Hopf algebras},  $K$-Theory  32  (2004),  no. 3, 231--251
%
\bibitem{Menichi2} Menichi, L., \textit{String topology for spheres},
arXiv:math/0609304.
%
\bibitem{Merkulov} Merkulov, S. A., \textit{De Rham model for string
topology}, Int. Math. Res. Not.  2004,  no. 55, 2955--2981.
%
\bibitem{Ng} Ng, L., \textit{Conormal bundles, contact homology, and knot invariants}, The interaction of finite type and GromovÐWitten invariants at the Banff International Research Station (2003), Geom. Topol. Monogr., vol. 8 (2006), pp. 129-144, arXiv:math.GT/0412330.
%
\bibitem{NouriThesis}
Nouri, Fereydoun, \textit{Even and odd graph homology (the commutative case)}, PhD dissertation, City Univ. of New York, 2004.
%
\bibitem{Polyakov} Polyakov, A. M., {\it Gauge Fields and Strings}, Chur: Harwood Academic Publishers, 1987. 
%
\bibitem{Quillen}
Quillen, D., \textit{Rational homotopy theory}, Ann. of Math., 90 (1969), 205-295.
%
\bibitem{SalamonWeber} Salamon, D. A., and Weber, J. \textit{Floer homology and the heat flow}, Geom. Funct. Anal. 16 (2006), no. 5, 1050--1138.
%
\bibitem{Schwarz} Schwarz, M., \textit{Morse homology}, Progress in Mathematics, vol. 111, BirkhŠuser Verlag, Basel, 1993.
%
\bibitem{Sta}
Stallings, J., {\it How not to prove the Poincar\'e conjecture.} Topology Seminar, Wisconsin, 1965, 
Ann. of Math. Studies, No. 60, Princeton, (1966). 
%
\bibitem{Stasheff}
Stasheff, J. D., \textit{Homotopy associativity of $H$-spaces. I, II}, Trans. Amer. Math. Soc. 108 (1963), 275-292; ibid. 108 1963 293-312. 
%
\bibitem{Su}
Sullivan, D., {\it Infinitesimal computations in topology}, Inst. Hautes ƒtudes Sci. Publ. Math. 47 (1977), 269-331. 
%
\bibitem{Su2}
Sullivan, D., {\it Geometric Topology: Localization, Periodicity and Galois Symmetry: The 1970 MIT Notes} K-Monographs in Mathematics, Dordecht: Springer, 2005.
%
\bibitem{Su3}
Sullivan, D., {\it Open and closed string field theory interpreted in classical algebraic topology} in {\it Topology, Geometry and Quantum Field Theory: Proceedings of the 2002 Oxford Symposium in Honour of the 60th Birthday of Graeme Segal} London Mathematical Society Lecture Note Series, Cambridge: Cambridge University Press, 2004.
%
\bibitem{Su4}
Sullivan, D., \textit{Sigma models and string topology} in \textit{Graphs And Patterns In Mathematics And Theoretical Physics: Proceedings Of The Stony Brook Conference On Graphs And Patterns In Mathematics And Theoretical Physics, Dedicated to Dennis Sullivan's 60th Birthday} Proceedings of Symposia in Pure Mathematics, American Mathematical Society, 2005.
%
\bibitem{SullivanSullivan}
Sullivan, D., and Sullivan, M., \textit{Open string topology and classical knots}, in progress.
%
%
\bibitem{Tamarkin}
Tamarkin, D., \textit{Another proof of M. Kontsevich formality theorem}, arXiv:math.QA/9803025. 
%
\bibitem{Thom}
 Thom, R. \textit{Les classes caractŽristiques de Pontrjagin des variŽtŽs triangulŽes},  Symposium internacional de topolog'a algebraica International symposium on algebraic topology 54-67 Universidad Nacional Aut—noma de MŽxico and UNESCO, Mexico City, 1958.
 %
 \bibitem{Toen}
 To\"en, B., \textit{The homotopy theory of $dg$-categories and derived Morita theory}, Invent. Math. 167 (2007), no. 3, 615-667.
 %
 \bibitem{TradlerThesis} Tradler, T. \textit{Poincar\'e Duality Induces a BV-Structure on Hochschild 
Cohomology}, PhD dissertation, City Univ. of New York, 2002.
%
\bibitem{Tradler} Tradler, T., \textit{The BV algebra on Hochschild
cohomology induced by infinity inner products}, arXiv:math.QA/0210150.
 %
 \bibitem{TZ} Tradler T. and Zeinalian, M., \textit{On the cyclic Deligne
conjecture}, J. Pure and Applied Algebra, 204 (2006), no. 2, 280-299.
%
\bibitem{Turaev} Turaev, V.,
{\it Skein quantization of Poisson algebras of
loops on surfaces}, Ann. Sci. Ecole Norm. Sup.
(4) 24 (1991), no. 6, 635-704.
%
\bibitem{Turaev2} Turaev, V.,
\textit{Algebras of loops on surfaces, algebras of knots, and quantization},
Braid group, knot theory and statistical mechanics, Adv. Ser. Math. Phys. 9 (1989), 59-95.
%
\bibitem{Vaintrob} Vaintrob, D., \textit{The string topology BV algebra,
Hochschild cohomology and the Goldman bracket on surfaces},
arXiv:math/0702859v1.
%
\bibitem{Vallette}
Valette, B., \textit{A Koszul duality for props}, arXiv:math/0411542
%
\bibitem{Viterbo} Viterbo, C., \textit{Functors and computations in Floer homology with applications. I}, Geom. Funct. Anal. 9 (1999), no. 5, 985-1033.
%
\bibitem{Viterbo'} Viterbo, C., \textit{Functors and computations in Floer homology with applications. II}, preprint.
%
\bibitem{Viterbo2} Viterbo, C., \textit{Exact Lagrange submanifolds, periodic orbits and the cohomology of free loop spaces}, J. Differential Geom. 47 (1997), no. 3, 420-468.
%
\bibitem{Viterbo3} Viterbo, C. and Hofer, H., \textit{The Weinstein conjecture in the presence of holomorphic spheres}, Comm. Pure Appl. Math. 45 (1992), no. 5, 583Ð622.
%
\bibitem{Wiesbrock}
Wiesbrock, H-W., \textit{The construction of the sh-Lie-algebra of closed bosonic strings}, 
Comm. Math. Phys. 145 (1992), no. 1, 17-42. 
%
\bibitem{WilsonThesis}
Wilson, S., \textit{On the Algebra and Geometry of a Manifold's Chains and Cochains}, PhD dissertation, State Univ. of New York, Stony Brook, 2005.
%
\bibitem{Wolpert} Wolpert, S.,
{\it On the Symplectic Geometry of Deformations
of Hyperbolic Surfaces} Ann. Math. 117
(1983), 207-234. 
%
\bibitem{Yang} Yang, T., \textit{A Batalin-Vilkovisky Algebra structure on
the Hochschild Cohomology of Truncated Polynomials}, arXiv:0707.4213.
%
\bibitem{Zapponi}
Zapponi, L., \textit{What is...a dessin d'enfant?} Notices of the Amer. Math. Soc. 50 (2003), no. 7, 788-789.
%
\bibitem{ZeinalianThesis}
Zeinalian, M., \textit{On some local combinatorial invariants of homology manifolds}
PhD dissertation, City Univ. of New York, 2001.
%

\end{thebibliography}
\end{document}